\documentclass[11pt,leqno,a4paper]{article}
\usepackage{mathrsfs}
\usepackage{amsmath}
\usepackage{amsfonts}
\usepackage{amsthm}
\usepackage{amssymb}
\usepackage[utf8]{inputenc}
\usepackage[T1]{fontenc}

\usepackage{xcolor, colortbl}
\usepackage{array,color}
\usepackage{arydshln}
\usepackage{threeparttable}
\usepackage{multirow}
\usepackage{float}
\usepackage{bm}
\usepackage{multirow}
\usepackage{bigdelim}

\usepackage[ruled,vlined]{algorithm2e}
\usepackage{algorithmic}
\usepackage{hhline}

\usepackage{upgreek}
\usepackage{tikz}
\usetikzlibrary{shadows}
\usetikzlibrary{shapes}

\newtheorem{theorem}{Theorem}[section]
\newtheorem{proposition}[theorem]{Proposition}

\theoremstyle{remark}
\newtheorem{remark}[theorem]{Remark}
\theoremstyle{definition}

\usepackage{authblk}

\title{Domain decomposition algorithms for two dimensional linear Schr{\"o}dinger equation}

\author[1]{Christophe Besse\thanks{christophe.besse@math.univ-toulouse.fr}}
\author[2]{Feng Xing\thanks{feng.xing@unice.fr}}

\affil[1]{Institut de Math{\'e}matiques de Toulouse UMR5219,
  Universit\'e de Toulouse; CNRS,
  UPS IMT, F-31062 Toulouse Cedex 9, France.}
\affil[2]{Maison de la Simulation, CEA Saclay France \& Laboratoire Paul Painlev{\'e}, Universit{\'e} Lille Nord de France.\\
  \emph{Present address:} Laboratoire J.A. Dieudonn{\'e}, Université de Nice \& Inria Sophia Antipolis, France.}

\begin{document}

\maketitle

\begin{abstract}
  This paper deals with two domain decomposition methods for two
  dimensional linear Schr\"{o}dinger equation, the Schwarz waveform
  relaxation method and the domain decomposition in space
  method. After presenting the classical algorithms, we propose a new
  algorithm for the Schr\"{o}dinger equation with constant potential
  and a preconditioned algorithm for the general Schr\"{o}dinger
  equation. These algorithms are studied numerically. The experiments show that
  the two new algorithms improve the  convergence rate and reduce
  the computation time. Besides the traditional Robin transmission
  condition, we also propose to use a newly constructed absorbing
  condition as the transmission condition.   
\end{abstract}



\section{Introduction}

The aim of this paper is to apply domain decomposition algorithms to the two dimensional linear Schr\"{o}dinger equation defined on $(0,T) \times \Omega$ with a real potential $V(t,x,y)$ 
\begin{equation} 
  \label{Schrodingereq}
  \left\{\begin{array}{ll}
      \mathscr{L}u := (i\partial_t  + \Delta + V) u = 0, & \ (t,x,y)\in (0,T)\times \Omega, \\
      u(0,x,y) = u_0(x,y), & \ (x,y) \in \Omega,
    \end{array} \right.
\end{equation}
where $\Omega = (x_l,x_r)\times (y_b,y_u)$ is a bounded spatial domain of $\mathbb{R}^2$ with $x_l,x_r,y_b,y_u \in \mathbb{R}$  and the initial datum $u_0 \in L^2(\Omega)$. 
The equation is complemented with homogeneous Neumann boundary condition on
bottom and top boundaries and Fourier-Robin boundary conditions
in left and right boundaries. They read:
\begin{align*}
  \partial_{\mathbf{n}} u = 0,\ y=y_b,y_u,\quad \partial_{\mathbf{n}} u + S_b u = 0,\  x=x_l,x_r,
\end{align*}
where $\partial_{\mathbf{n}}$ denotes the normal derivative, $\mathbf{n}$ being the outwardly unit vector on the boundary $\partial \Omega$, and the operator $S_b$ is some  transmission operator.

We consider in this paper two domain decomposition methods. The
first one is the Schwarz waveform relaxation method without
overlap (SWR) \cite{Martin2005,Gander2003_wave}, which is based
on the time-space domain decomposition. The time-space domain
$(0,T)\times \Omega$ is decomposed into some subdomains
$(0,T)\times \Omega_j$, $j=1,2,...,N$. The solution is computed
on each subdomain and the time-space boundary values are
transmitted via transmission conditions. The derivation of
efficient transmission conditions is one of the key points of the
SWR method. For Schr\"{o}dinger equation, some transmission
conditions are proposed in
\cite{Halpern2010_sch,Antoine2014,XF20151d}, such as Robin
transmission condition, optimal transmission condition etc.. 

The second method we consider here is the domain decomposition in space method (DDS)
\cite{Cai1994,Wu1998}. First of all, the time dependent equation is
semi-discretized in time with an implicit scheme on the entire spatial
domain. This procedure leads to a stationary equation in space. Then
standard domain decomposition methods (such as the optimized Schwarz
method \cite{Gander2006osw,Boubendir2012,NATAF1995}) are applied
to this stationary equation. The DDS method demands a conforming
time discretization. The use of nonconforming discretization in time is
non standard for the Schr\"{o}dinger equation and we therefore
fulfill this requirement.

We propose in this article to show the effectiveness of new
transmissions conditions (expressed in term of absorbing transparent
conditions) when we apply the two classical algorithms to the
  Schr\"odinger equation. The study of the interface problem allows us
  to introduce some new algorithms which
  significantly reduce both
  the computational time and the number of iterations. We compare them
  to the classical widely used Robin transmission
  condition with various intensive numerical tests made on parallel
  computers with  up to 1024 subdomains.

This paper is organized as follows. In Section \ref{Sec_DDAlgo},
we present the classical SWR and DDS
algorithms. We show how the classical DDS algorithm can be interpreted as a
combination of some classical SWR algorithms. In Section
\ref{Sec_discreteInterfaceproblem}, we construct an interface
problem and analyse its properties. The discretization of the Schr\"odinger equation is also
provided. Based on these properties, we propose new algorithms
and preconditioned algorithms in the two following sections. In
Section \ref{Sec_NumResults}, we provide numerical experiments
which show the efficiency of our new algorithms. A conclusion is drawn in the
last section. 

\section{Domain decomposition algorithms}
\label{Sec_DDAlgo}

\subsection{Geometric configuration} 
The interval $(x_l, x_r)$ is divided into $N$ subintervals $(a_j,
b_j)$ without overlap. The points $a_j$ and $b_j$ denote the ends
of the subintervals $(a_j, b_j)$. Thus, the entire domain
$\Omega$ is decomposed into $N$ non overlapping subdomains
$\Omega_j = (a_j, b_j) \times (y_b, y_u)$, $j = 1, 2,..., N$ (see
Figure \ref{fig_omega0} for $N = 3$). We denote the normal
derivative on subdomain $\Omega_j$ by $\partial_{\mathbf{n}_j}$. 

\definecolor{col1}{rgb}{0.8,1,0.8}
\definecolor{col2}{rgb}{1,0.8,0.8}
\definecolor{col3}{rgb}{0.8,0.8,1}
\begin{figure}[!htbp]
  \centering
  \begin{tikzpicture}
    
    \draw [color=gray,fill=col1] (0,0) rectangle (2.6,2);
    \draw [color=gray,fill=col2] (2.6,0) rectangle (5.2,2);
    \draw [color=gray,fill=col3] (5.2,0) rectangle (7.8,2);
    \draw[>=stealth,->] (-0.5,1) -- (8.3,1);
    \draw (8.1,0.8) node[right] {$x$};
    \draw [>=stealth,->] (3.9,-0.5) -- (3.9,2.5);
    \draw (3.7,2.2) node[above] {$y$};    
    
    \draw[>=stealth,->] (2.6,1.8) -- (2.3,1.8);
    \draw (2.2,1.8) node[below,scale=0.7] {$\mathbf{n}_2$};

    \draw[>=stealth,->] (5.2,1.8) -- (5.5,1.8);
    \draw (5.6,1.8) node[below,scale=0.7] {$\mathbf{n}_2$};    
    
    \draw (0.2,-0.1) node[below] {$x_l=a_1$};
    \draw (2.6,0) node[below] {$b_1=a_2$};
    \draw (5.2,0) node[below] {$b_2=a_3$};
    \draw (7.8,0) node[below] {$b_2=x_r$};
    \draw (1.0,1.3) node[scale=0.8] {$\Omega_1$};
    \draw (3.6,1.3) node[scale=0.8] {$\Omega_2$};
    \draw (6.1,1.3) node[scale=0.8] {$\Omega_3$};
  \end{tikzpicture}
  \caption{Geometric configuration.}
  \label{fig_omega0}
\end{figure}
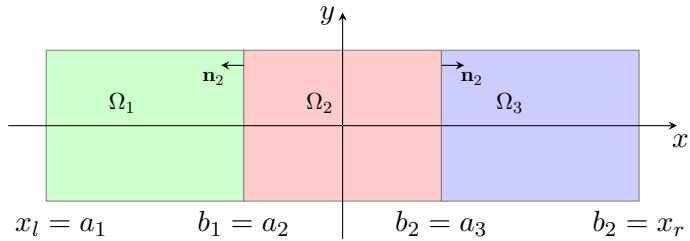

\begin{figure}[!htbp]
  \centering
  \includegraphics[width=0.35\textwidth]{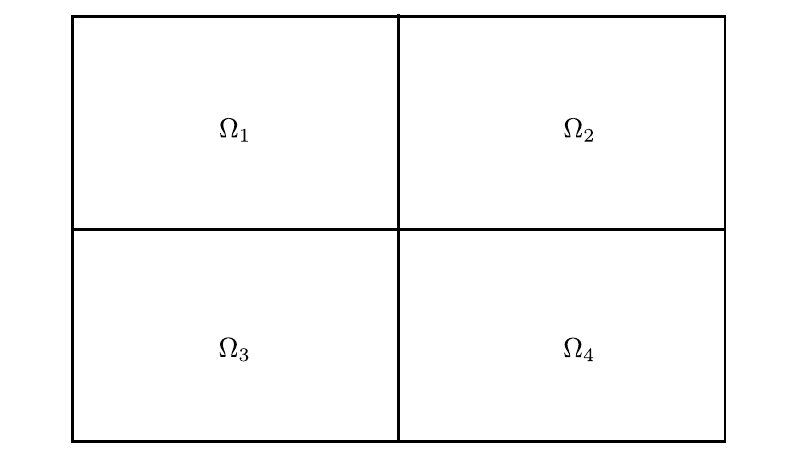}
  \includegraphics[width=0.35\textwidth]{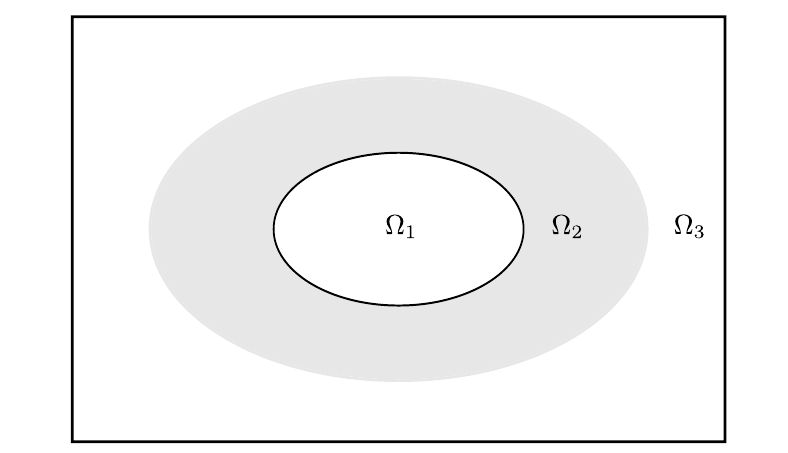}
  \caption{Two other ways of domain decomposition.}
  \label{fig_ddespace}
\end{figure}

There are obviously other ways to decompose the entire
domain. One way is illustrated in Figure \ref{fig_ddespace}
(left) for $N = 4$. The intervals $(x_l, x_r)$ and $(y_b, y_u)$
are simultaneously decomposed into subintervals in both spatial
directions. In this configuration, an artificial cross point
appears. It is well known that the domain decomposition method
with cross points is a difficult problem since the
problem becomes singular at this point. Another possibility is
illustrated in Figure \ref{fig_ddespace} (right) for $N = 3$. The
entire domain is decomposed into an ellipsis and some rings. This
approach has many disadvantages for parallel computing. Indeed,
we would like to control the number of cells for the meshes of
each subdomains. Their sizes have to be equivalent to insure a
good balance between different process. Thus, we restrict ourselves in this paper
to the first description (see Figure \ref{fig_omega0}). 


\subsection{Classical SWR algorithm}
\label{Sec_dist}
The Schwarz methods being iterative, we label the iteration number
of the algorithms by $k$ and denotes by $u_j^{k}$ the solution on subdomain $(0,T) \times \Omega_j$ at iteration $k = 1,2,...$.

The classical SWR algorithm is given by
\begin{equation} 
  \label{AlgoSWR}
  \left\{ 
    \begin{array}{ll}
      \mathscr{L} u_j^{k} = 0, & (t,x,y) \in (0,T) \times \Omega_j, \  \\
      u^{k}_j(0,x,y) = u_0(x,y), & (x,y) \in \Omega_j, \\
      B_j u_j^{k} = B_j u_{j-1}^{k-1}, & x = a_j,\\
      B_j u_j^{k} = B_j u_{j+1}^{k-1}, & x = b_j,\\
      \partial_{\mathbf{n}_j} u_j^{k} = 0, & y=y_l,y_b,
    \end{array}  
  \right. 
\end{equation}
with a special treatment for the two extreme subdomains
$(0,T) \times \Omega_1$ and $(0,T) \times \Omega_N$ since
the boundary conditions are imposed on $x=a_1$ and $x=b_N$ 
\begin{displaymath}
  B_1 u_1^{k}=(\partial_{\mathbf{n}}+S_b)u_1^k = 0, \ x=a_1, \quad B_N u_N^{k}=(\partial_{\mathbf{n}}+S_b)u_N^k  = 0, \ x=b_N.
\end{displaymath}
The boundary condition at interface nodes $a_j$ and $b_j$ is given in term of operator $B_j$ defined by
\begin{equation}
  \label{condBj} 
  B_j = \partial_{\mathbf{n}_j} + S_j,\ j=1,...,N,
\end{equation}
where $S_j$ is a transmission operator. Various $S_j$
operators can be considered. The 
classical widely used Robin transmission condition  is
given by 
\begin{equation}
  \label{Cond_Robin}
  S_j =  - ip,\ p \in \mathbb{R}^{+}, \ j=1,2,...,N.
\end{equation}
Traditionally, the optimal transmission operator is
  given in term of transparent boundary conditions. For the
  general linear two dimensional Schr\"odinger equation, we only have access to
  approximated version of the TBCs given by the recently
  constructed absorbing boundary condition
  $S_{\mathrm{pade}}^m$ \cite{Antoine2012,Antoine2013}
  which we used as the transmission condition
\begin{equation}
  \label{Cond_S2p}
  S_j = -i \sqrt{i \partial_t + \Delta_{\Gamma_j} + V}, \ j=1,2,...,N,
\end{equation}
where $\Gamma_1=\{b_1\} \times (y_b,y_u)$,
$\Gamma_j=\{a_j,b_j\}  \times (y_b,y_u), j=2,3,...,N-1$ and
$\Gamma_N=\{a_N\} \times (y_b,y_u)$. In our case, the
Laplace–Beltrami operator $\Delta_{\Gamma_j}$ is
$\partial_y^2$. 
Numerically, this operator is approximated by Pad\'{e} approximation of order $m$
\begin{displaymath}
  \sqrt{i\partial_t + \Delta_{\Gamma_j} + V} u \approx
  \Big( \sum_{s=0}^m a_s^m - \sum_{s=1}^m a_s^m d_s^m (i\partial_t + \Delta_{\Gamma_j} + V + d_s^m)^{-1} \Big) u.
\end{displaymath}
If the potential $V=0$ and the spatial domain is
$\mathbb{R}^2$, the absorbing boundary condition is
actually a transparent boundary conditions. Then it could be proven that the transmission
condition $S_{\mathrm{pade}}^m$ leads to an optimal SWR
method. 

The classical SWR algorithm is initialized by an initial
  guess of $B_j u_j^0 |_{x=a_j,b_j}$, $j=1,2,...,N$. The
  boundary conditions for any subdomain $\Omega_j$ at
  iteration $k+1$ involve the knowledge of the values of
  the functions on adjacent subdomains $\Omega_{j-1}$ and
  $\Omega_{j+1}$ at prior iteration $k$. Thanks to the
  initial guess, we can \emph{solve} the Schr\"odinger equation on
  each subdomain, allowing to build the new boundary
  conditions for the next step, \emph{communicating} them to other
  subdomains. This procedure is summarized 
  in \eqref{AlgoSWRGraph} for $N=3$ subdomains at iteration
  $k$.
\begin{equation}
  \label{AlgoSWRGraph}
  \begin{pmatrix}
    B_1 u_1^{k} |_{x=b_1}\\
    B_2 u_2^{k} |_{x=a_2}\\
    B_2 u_2^{k} |_{x=b_2}\\
    B_3 u_3^{k} |_{x=a_3}
  \end{pmatrix}
  \underrightarrow{\hspace{0.1cm} \text{\footnotesize{\emph{Solve}}} \hspace{0.1cm}}
  \begin{pmatrix}
    u_1^{k} \\
    u_2^{k} \\
    u_3^{k}
  \end{pmatrix}
  \rightarrow
  \begin{pmatrix}
    B_2 u_1^{k} |_{x=b_1}\\
    B_1 u_2^{k} |_{x=a_2}\\
    B_3 u_2^{k} |_{x=b_2}\\
    B_2 u_3^{k} |_{x=a_3}
  \end{pmatrix}
  \underrightarrow{ \hspace{0.1cm} \text{\footnotesize{\emph{Comm.}}} \hspace{0.1cm} }
  \begin{pmatrix}
    B_1 u_1^{k+1} |_{x=b_1}\\
    B_2 u_2^{k+1} |_{x=a_2}\\
    B_2 u_2^{k+1} |_{x=b_2}\\
    B_3 u_3^{k+1} |_{x=a_3}
  \end{pmatrix}.
\end{equation}
Let us define the flux at iteration $k$ by 
$$g^k = (B_1 u_1^{k} |_{x=b_1}, \cdots, B_j u_j^{k} |_{x=a_j}, B_j u_j^{k} |_{x=b_j}, \cdots, B_N u_N^{k} |_{x=a_N})^{\top},$$
"$\cdot^\top$" denoting the transpose of the vector. Thanks
  to this definition, we give a new interpretation to the
  algorithm which can be written as
\begin{equation}
  \label{InterfaceAlgoContinues}
  g^{k+1} = \mathcal{R}_{c} g^{k},
\end{equation}
where $\mathcal{R}_c$ is a linear operator. The solution to
  this iteration process is given as the solution to the continuous interface problem
\begin{equation}
  \label{InterfacePbContinues}
  (I-\mathcal{R}_{c}) g = 0,  \ g=\lim_{k \rightarrow \infty} g^k.
\end{equation}
where  $I$ is identity operator.  

\subsection{Classical DDS algorithm \label{sec:dds}}
\label{Sec_dds}

The other algorithm we consider in this paper is the domain
decomposition in space algorithm (DDS). The equation
\eqref{Schrodingereq} is first semi-discretized in time on the
entire domain $(0,T) \times \Omega$. The time interval
$(0,T)$ is discretized uniformly with $N_T$ intervals
of length $\Delta t$. We denotes $u_{n}$ (resp. $V_n$) an
  approximation of the solution $u$ (resp. $V$) at time $t_n=n\Delta t$.
The Crank-Nicolson scheme on $(0,T) \times \Omega$ reads 
\begin{displaymath}
  i\frac{u_{n} - u_{n-1}}{\Delta t} + \Delta \frac{u_{n} +
    u_{n-1}}{2} + \frac{V_{n} + V_{n-1}}{2} \frac{u_{n} +
    u_{n-1}}{2} = 0, \ 1 \leqslant n \leqslant N_T. 
\end{displaymath} 
By introducing new variables $v_{n} = (u_{n} + u_{n-1})/2$
with $v_{0} = u_{0}$ and $W_{n} =  (V_{n} + V_{n-1})/2$, we
get a stationary equation defined on $\Omega$ with unknown
$v_n$ 
\begin{equation}
  \label{EqSemi}
  \mathcal{L}_{\mathbf{x}} v_{n} = \frac{2i}{\Delta t} u_{n-1},
\end{equation}
where $\mathcal{L}_{\mathbf{x}}:= \frac{2i}{\Delta t} + \Delta  + W_n$. 
  We recover the original unknown by $u_{n} = 2 v_{n} - u_{n-1}$.
%
%
Then, the optimized Schwarz algorithm is applied to the
stationary equation \eqref{EqSemi}. We denote by
$R_j,j=1,2,...,N$ the restriction operator from $\Omega$ to
$\Omega_j$. At time $t_n$, the classical algorithm reads 
\begin{equation} 
  \label{AlgoDDS}
  \left\{ \begin{array}{ll}
      \mathscr{L}_{\mathbf{x}} v_{n,j}^{k} = \frac{2i}{\Delta t} R_j u_{n-1}, & (x,y) \in  \Omega_j,\\
      \partial_{\mathbf{n}_j} v_{n,j}^{k} + \overline{S}_j v_{n,j}^{k} = \partial_{\mathbf{n}_{j}} v_{n,j-1}^{k-1} +  \overline{S}_j v_{n,j-1}^{k-1}, & x = a_j,\\
      \partial_{\mathbf{n}_j} v_{n,j}^{k} + \overline{S}_j v_{n,j}^{k} = \partial_{\mathbf{n}_j} v_{n,j+1}^{k-1} + \overline{S}_j v_{n,j+1}^{k-1}, & x = b_j,\\
      \partial_{\mathbf{n}_j} v_{n,j}^{k} = 0, & y=y_l,y_b,
    \end{array}  \right. 
\end{equation}
where $v_{n,j}^k$ denotes the unknown at time $t_n$, on
  subdomain $\Omega_j$ at iteration $k$ and $\overline{S}_j$ is the semi-discrete form of
  $S_j$ given by \eqref{Cond_Robin} or \eqref{Cond_S2p}. 
A special treatment for the two extreme subdomains is
needed and the boundary conditions read
$\partial_{\mathbf{n}_1} v_{n,1}^{k} + \overline{S}_1
v_{n,1}^{k} = 0,x=a_1$, $\partial_{\mathbf{n}_N}
v_{n,N}^{k} + \overline{S}_N v_{n,N}^{k} =
0,x=b_N$. 

Since the interval $(t_{n-1},t_n)$ contains only one
time step, the DDS algorithm can be numerically
interpreted as a sequence of some SWR algorithms.
\begin{center}
  \begin{minipage}{0.85\textwidth}
    \centering
    \begin{algorithm}[H]
      \label{AlgoDDS_swr}
      \caption{DDS algorithm}
      The initial datum is $u_{0}$. \\
      \For{$n=1,2,...,N_T$}{
        \vbox{}
        Apply the SWR algorithm to
        \begin{displaymath}
          \left \{ \begin{array}{ll}
              \mathscr{L} u = 0, \ (t,x,y) \in (t_{n-1},t_n) \times \Omega,\\
              u(0,x,y) = u_{n-1}(x,y), \ (x,y) \in \Omega, \\
            \end{array} \right. 
        \end{displaymath} 
        where $t_n = n \Delta t$.        
      }	 
    \end{algorithm}
  \end{minipage}
\end{center}

\section{Discrete interface problem}
\label{Sec_discreteInterfaceproblem}

We know by \eqref{InterfacePbContinues} that the classical SWR algorithm reduces to the interface problem
  \begin{equation*}
    (I-\mathcal{R}_{c}) g = 0,  \ g=\lim_{k \rightarrow \infty} g^k.
  \end{equation*}
  The aim of this section is discretize this relation and to show that the discrete interface problem can be written as
\begin{equation}
  \label{InterfacePbdisctILgd}
  (I - \mathcal{L}_h) \mathbf{g} = \mathbf{d},
\end{equation}
where the vector $\mathbf{g}$ is the discrete form of $g$, $\mathbf{d}$ is a vector
\begin{displaymath}
  \mathbf{d} = (\mathbf{d}_{1,r}^{\top}, \mathbf{d}_{2,l}^{\top}, \mathbf{d}_{2,r}^{\top},\cdots, \mathbf{d}_{N,l}^{\top})^{\top} \in \mathbb{R}^{N_y \times N_T},
\end{displaymath}       
and $\mathcal{L}_h$ is a block matrix (the notation ``MPI $j$''
  above the columns of the matrix will be
used in section \ref{new_algo_sec})
\begin{equation}
  \label{matL}
  \mathcal{L}_{h} = 
  \begin{pmatrix}
    \multicolumn{1}{l}{\overbrace{\hspace{2.0em}}^{\mathrm{MPI}\ 0}} & 
    \multicolumn{2}{l}{\overbrace{\hspace{5.0em}}^{\mathrm{MPI}\ 1}} &
    \multicolumn{2}{l}{\overbrace{\hspace{5.0em}}^{\mathrm{MPI}\ 2}} &
    & 
    \multicolumn{2}{l}{\overbrace{\hspace{8.0em}}^{\mathrm{MPI}\ N-2}} &
    \multicolumn{2}{l}{\overbrace{\hspace{2.0em}}^{\mathrm{MPI}\ N-1}}
    \\
    & X^{2,1} & X^{2,2} & & & \\
    X^{1,4} \\
    & & & X^{3,1} & X^{3,2} \\
    & X^{2,3} & X^{2,4} \\
    & & & & & \cdots \\
    & & & X^{3,3} & X^{3,4} \\
    & & & & & & X^{N-1,1} & X^{N-1,2}\\
    & & & & &\cdots \\
    & & & & & & & & X^{N,1}\\
    & & & & & & X^{N-1,3} & X^{N-1,4}
  \end{pmatrix}.
\end{equation}                                                     

The concrete definitions of $\mathbf{d}_{j,l}$, $\mathbf{d}_{j,r}$,
$j=1,2,...,N$ and the bloks in $\mathcal{L}_h$ are given in
proposition \ref{InterfGlobalRobin} and  proposition
\ref{InterfGlobalS2P}. The subsection \ref{prop_Lh} is devoted to the
  derivation of some important properties of $\mathcal{L}_h$, which give
  the ground to make the new algorithm given in section \ref{new_algo_sec}.


\subsection{Preliminaries related to discretization}
\label{Sec_discrete}
Without loss of generality, we present the discretization of \eqref{Schrodingereq} on $(0,T) \times \Omega$ with the following boundary conditions \cite{Antoine2012,Antoine2013}
\begin{equation*} 
  \left\{ \begin{array}{ll}
      \partial_{\mathbf{n}} u = 0,\ y=y_b,y_u, \\ 
      \partial_{\mathbf{n}} u + S u = l(t,y),\  x=x_l, \\
      \partial_{\mathbf{n}} u + S u = r(t,y),\  x=x_r,
    \end{array}  \right. 
\end{equation*}
where $l(t,y)$ and $r(t,y)$ are two functions. The
  discrete version of \eqref{AlgoSWR} follows
  immediately. Concerning the semi-discretization with
  respect to time, we follow the procedure given in
  section \ref{sec:dds} for the classical DDS
  algorithm. The scheme therefore reads
\begin{equation}
  \frac{2i}{\Delta t} v_{n} + \Delta v_{n} + W_n v_{n} = \frac{2i}{\Delta t} u_{n-1}.
\end{equation}
The semi-discrete transmission condition is given by
\begin{displaymath}
  \partial_{\mathbf{n}} v_{n} + \overline{S} v_{n} = l_{n}, \ x=x_l, \quad
  \partial_{\mathbf{n}} v_{n} + \overline{S} v_{n} = r_{n}, \ x=x_r,
\end{displaymath}
where $l_n=l(n\Delta t,y)$, $r_n=r(n\Delta t,y)$ and
$\overline{S}$ is the semi-discrete form of
$S$. If we consider the Robin transmission
  condition \eqref{Cond_Robin}, we have
  \begin{equation}
    \mathrm{Robin: } \quad \overline{S} v_{n} = -ip
    \cdot v_{n}. \label{CondSemi_Robin} 
  \end{equation}
  The approximation of the transmission condition
  \eqref{Cond_S2p} is given by
  \begin{equation}
    S_{\mathrm{pade}}^m: \quad \overline{S} v_{n} = -i \displaystyle \sum_{s=0}^m a_s^m v_{n} + i\sum_{s=1}^m a_s^m d_s^m \varphi^{n-1/2}_{s}, \label{CondSemi_S2p}
  \end{equation}
  where $a_s^m = e^{i\theta/2}/(m\cos^2(\frac{(2s-1)\pi}{4m}))$,
  $d_s^m = e^{i\theta} \tan^2 (\frac{(2s-1)\pi}{4m})$, $s=0,2,...,m$,
  $\theta = \frac{\pi}{4}$. The auxiliary functions $\varphi^{n-1/2}_{s}$,
  $s=1,2,...,m$ are defined as the solutions of the set of
  equations 
\begin{displaymath}
  \left \{ \begin{aligned}
      & \big( \frac{2i}{\Delta t} + \Delta_{\Gamma} + W_{n} + d_s^m \big) \varphi^{n-1/2}_{s} - v_{n} = \frac{2i}{\Delta t} \varphi^{n-1}_{s},  \\
      & \varphi_{s}^{n} = 2 \varphi_{s}^{n-1/2} - \varphi_{s}^{n-1}, \ \varphi_{s}^0 = 0.
    \end{aligned} \right. 
\end{displaymath}

The spatial approximation is realized by the standard $\mathbb{Q}_1$
finite element method. The uniform mesh size of a discrete element is
$(\Delta x, \Delta y)$. We denote by $N_x$ (resp. $N_y$) the number of
nodes in $x$ (resp. $y$) direction on each subdomain. Let us denote by
$\mathbf{v}_{n}$ (resp. $\mathbf{u}_{n}$) the nodal interpolation
vector of $v_{n}$ (resp. $u_{n}$), $\mathbf{l}_{n}$
(resp. $\mathbf{r}_{n}$) the nodal interpolation vector of $l_{n}$
(resp. $r_{n}$), $\mathbb{M}$ the mass matrix, $\mathbb{S}$ the
stiffness matrix and $\mathbb{M}_{W_n}$ the generalized mass matrix
with respect to $\int_{\Omega} W_n v \phi dx$. Let
$\mathbb{M}^{\Gamma}$ the boundary mass matrix, $\mathbb{S}^{\Gamma}$
the boundary stiffness matrix and $\mathbb{M}^{\Gamma}_{W_{n}}$ the
generalized boundary mass matrix with respect to $\int_{\Gamma} W_n v
\phi d\Gamma$. We denote by $Q_{l}$ (resp. $Q_{r}$) the restriction
operators (matrix) from $\Omega$ to $\{x_l\} \times (y_b,y_u)$
(resp. $\{x_r\} \times (y_b,y_u)$) and $Q^{\top}=(Q_{l}^{\top},
Q_{r}^{\top})$. The matrix formulation for the transmission condition
Robin is therefore given by 
\begin{equation}
  \label{NproblemRobin}
  \mathrm{Robin }: \quad \Big( \mathbb{A}_{n} + ip \cdot \mathbb{M}^{\Gamma} \Big)
  \mathbf{v}_{n} =  \frac{2i}{\Delta t} \mathbb{M}
  \mathbf{u}_{n-1} - \mathbb{M}^{\Gamma}  Q^{\top} 
  \begin{pmatrix}
    \mathbf{l}_{n} \\ 
    \mathbf{r}_{n}
  \end{pmatrix},   
\end{equation}
where $\mathbb{A}_{n} = \frac{2i}{\Delta t} \mathbb{M} - \mathbb{S} + \mathbb{M}_{W_{n}}$. The size of this linear system is $N_x \times N_y$. If we consider the transmission condition $S_{\mathrm{pade}}^m$, we have
\begin{gather}
  \begin{pmatrix}
    \mathbb{A}_{n} + i(\sum_{s=0}^m a_s^m)  \cdot
    \mathbb{M}^{\Gamma} & \mathbb{B}_{1} & \mathbb{B}_{2} & \cdots & \mathbb{B}_{m} \\
    \mathbb{C} & \mathbb{D}_{1}^n  \\
    \mathbb{C} & & \mathbb{D}_{2}^n  \\
    \vdots & & & \ddots \\
    \mathbb{C} & &  & & \mathbb{D}_{m}^n 
  \end{pmatrix}
  \begin{pmatrix}
    \mathbf{v}_{n}\\
    \bm{\varphi}^{n-1/2}_{1}\\
    \bm{\varphi}^{n-1/2}_{2}\\
    \vdots\\
    \bm{\varphi}^{n-1/2}_{m}
  \end{pmatrix} \nonumber
  \\
  = \begin{pmatrix}
    \mathbb{M}_j \\
    & Q \mathbb{M}^{\Gamma} Q^{\top} \\
    & & Q \mathbb{M}^{\Gamma} Q^{\top} \\
    & & & \ddots \\
    & & & & Q \mathbb{M}^{\Gamma} Q^{\top}
  \end{pmatrix}
  \begin{pmatrix}
    \mathbf{u}_{n-1}\\
    \bm{\varphi}^{n-1}_{1}\\
    \bm{\varphi}^{n-1}_{2}\\
    \vdots\\
    \bm{\varphi}^{n-1}_{m}
  \end{pmatrix}
  - 
  \begin{pmatrix}
    \mathbb{M}^{\Gamma} Q^{\top}
    \begin{pmatrix}
      \mathbf{l}_{n}\\
      \mathbf{r}_{n}
    \end{pmatrix}\\
    0\\
    \vdots \\
    0
  \end{pmatrix},
  \label{NproblemSP}
\end{gather}  
with
\begin{align*}
  & \mathbb{B}_{s} = -i a_s^m d_s^m \mathbb{M}^{\Gamma} Q^{\top} , \ 1
  \leqslant s \leqslant m,\\
  & \mathbb{C} = - Q \mathbb{M}^{\Gamma}, \\
  & \mathbb{D}_{s}^n = Q (\frac{2i}{\Delta
    t} \mathbb{M}^{\Gamma} -  \mathbb{S}^{\Gamma} +
  \mathbb{M}^{\Gamma}_{W_{n}} + d_s^m \mathbb{M}^{\Gamma}) Q^{\top},  \ 1
  \leqslant s \leqslant m.
\end{align*}                 
%
%
It is a linear system with unknown $(\mathbf{v}_{n},
\bm{\varphi}^{n-1/2}_{1}, ..., \bm{\varphi}^{n-1/2}_{m})$ where
$\bm{\varphi}^{n-1/2}_{s}$ is the nodal interpolation of
$\varphi^{n-1/2}_{s}$ on the boundary. 

\begin{remark}
  The $S_{\mathrm{pade}}^m$ transmission condition involves a larger
    linear system to solve than the one of the Robin transmission
    condition. The cost of the algorithm with the
    $S_{\mathrm{pade}}^m$ transmission condition is therefore more expensive. 
\end{remark}

The equations \eqref{NproblemRobin} and \eqref{NproblemSP} are given
  for a fixed discrete time $t_n$. They however can be
  written globally in time by some straightforward calculations. 
\begin{proposition}
  For the Robin transmission condition, the global form in time of the equation \eqref{NproblemRobin} is
  \begin{equation}
    \label{NproblemGlobalRobin}
    (\mathbf{A}  - \mathbf{B}) \mathbf{v} =
    \mathbf{F} - \mathbf{M}^{\Gamma} \mathbf{Q}^{\top} \mathbf{g}, 
  \end{equation}
  where  $\mathbf{B} = -ip \cdot \mathbf{M}^{\Gamma} = -ip \cdot \mathrm{diag}_{N_T}\{ \mathbb{M}^{\Gamma} \}$, $\mathbf{Q}_{l}^{\top} = \mathrm{diag}_{N_T} \{Q_{l}^{\top}\}$, $\mathbf{Q}_{r}^{\top} = \mathrm{diag}_{N_T} \{Q_{r}^{\top}\}$ and
  \begin{gather*}
    \mathbf{A} = \begin{pmatrix}
      \mathbb{A}_{1} &  &  &  &\\
      -\frac{4i}{\Delta t}\mathbb{M} & \mathbb{A}_{2} &  &  & \\
      \frac{4i}{\Delta t} \mathbb{M} & - \frac{4i}{\Delta t} \mathbb{M} & \mathbb{A}_{3} &  &  & \\
      \vdots & \vdots & & \ddots & \\
      & & & - \frac{4i}{\Delta t} \mathbb{M} & \mathbb{A}_{N_T}
    \end{pmatrix},\
    \mathbf{v} =
    \begin{pmatrix}
      \mathbf{v}_{1}\\
      \mathbf{v}_{2}\\
      \vdots \\
      \vdots \\
      \mathbf{v}_{N_T}
    \end{pmatrix},\\
    \mathbf{F} = \frac{2i}{\Delta t}
    \begin{pmatrix}
      \mathbb{M} \mathbf{u}_{0} \\
      -\mathbb{M} \mathbf{u}_{0}\\
      \vdots \\
      (-1)^{N_T-1}\mathbb{M} \mathbf{u}_{0}
    \end{pmatrix},\    
    \mathbf{Q}^{\top} = \begin{pmatrix}
      Q_{l}^{\top} &  &  &  &  Q_{r}^{\top} \\
      & Q_{l}^{\top} &  &  &  &  Q_{r}^{\top} \\
      & & \ddots &  &  &  & \ddots \\
      & & & Q_{l}^{\top} &  &  &  &  Q_{r}^{\top}
    \end{pmatrix}.
  \end{gather*}
\end{proposition}

\begin{proposition}
  If we consider the transmission condition $S_{\mathrm{pade}}^m$,
  then the equation \eqref{NproblemSP} can be written globally in time
  as \eqref{NproblemGlobalRobin} 
  \begin{displaymath}
    (\mathbf{A}  - \mathbf{B}) \mathbf{v} =
    \mathbf{F} - \mathbf{M}^{\Gamma} \mathbf{Q}^{\top} \mathbf{g}, 
  \end{displaymath}  
  with $\mathbf{B}$ given by 
  \begin{equation}
    \label{BjSP}
    \renewcommand{\arraystretch}{1.2}
    \mathbf{B} = -
    \begin{pmatrix}
      c_a \mathbb{M}^{\Gamma} + \mathbb{Y}^{1,1} \\
      \mathbb{Y}^{2,1} & c_a \mathbb{M}^{\Gamma} + \mathbb{Y}^{2,2} \\
      \mathbb{Y}^{3,1} & \mathbb{Y}^{3,1} & c_a \mathbb{M}^{\Gamma} + \mathbb{Y}^{3,3} \\
      \vdots & \vdots & & \ddots  \\
      \mathbb{Y}^{N_T,1} & \mathbb{Y}^{N_T,2} & \mathbb{Y}^{N_T,3} & \cdots  & c_a \mathbb{M}^{\Gamma} + \mathbb{Y}^{N_T,N_T}
    \end{pmatrix},
  \end{equation}
  where $c_a=i(\sum_{s=0}^m a_s^m)$.
  
\end{proposition}

\proof
  According to \eqref{NproblemSP}, for $s=1,2,...,m$, $n=1,2,...,N_T$, we have  
    \begin{equation*}
      - Q \mathbb{M}^{\Gamma} \mathbf{v}_{n} +  \mathbb{D}_{s}^n \bm{\varphi}_{s}^{n-1/2} = \frac{2i}{\Delta t} Q_j \mathbb{M}^{\Gamma} Q^{\top} \bm{\varphi}_{s}^{n-1},\  \ n \geqslant 1,
    \end{equation*}
    where $\mathbb{D}_{s}^n = Q (\frac{2i}{\Delta t} \mathbb{M}^{\Gamma}
    -  \mathbb{S}^{\Gamma} + \mathbb{M}^{\Gamma}_{W_{n}} + d_s^m
    \mathbb{M}^{\Gamma}) Q^{\top}$. We therefore obtain
    \begin{equation*}
      \bm{\varphi}_{s}^{n-1/2} = (\mathbb{D}_{s}^n)^{-1} Q_j \mathbb{M}^{\Gamma} \mathbf{v}_{n} + ( \mathbb{D}_{s}^n)^{-1} \frac{2i}{\Delta t} Q \mathbb{M}^{\Gamma} Q^{\top} \bm{\varphi}_{s}^{n-1},
    \end{equation*}
    and
    \begin{equation*}
      \bm{\varphi}_{s}^{n-1} = 2\bm{\varphi}_{s}^{n-3/2} - \bm{\varphi}_{s}^{n-2},\ \bm{\varphi}_{s}^{0} = 0, \ n \geqslant 2.
    \end{equation*}
  %
  By induction, $\bm{\varphi}_{s}^{n-1/2}$ is given by
  \begin{equation}
    \label{phiLv}
    \bm{\varphi}_{s}^{n-1/2} = \sum_{p=1}^{n} \mathbb{L}_{s}^{n,p} \mathbf{v}_{p},
  \end{equation}
  where $\mathbb{L}_{s}^{n,p}$ are matrix. Replacing
  $\bm{\varphi}_{s}^{n-1/2}$ by $\mathbf{v}_{p}$ in the
  first row of the equation \eqref{NproblemSP} gives
  \begin{align}
    \label{AvulSP}
    \Big(  \mathbb{A}_{n}  + i(\sum_{s=0}^m a_s^m)  \cdot
    \mathbb{M}^{\Gamma} \Big)
    \mathbf{v}_{n}  & + \sum_{p=1}^{n}  \mathbb{Y}^{n,p} \mathbf{v}_{p}
    =  \frac{2i}{\Delta t} \mathbb{M}
    \mathbf{u}_{n-1}
    - \mathbb{M}^{\Gamma} Q^{\top}
    \begin{pmatrix}
      \mathbf{l}_{n}\\
      \mathbf{r}_{n}
    \end{pmatrix},
  \end{align}
  where $\mathbb{Y}^{n,p} := -i \sum_{s=1}^m a_s^m d_s^m
  \mathbb{M}^{\Gamma } Q^{\top} \mathbb{L}_{s}^{n,p}$ 
  Then, according to \eqref{AvulSP}, the matrix
    $\mathbf{B}$ is given by \eqref{BjSP}.

\begin{remark}
  Following the same procedure, the equation \eqref{AlgoSWR} is discretized on each
  subdomain $(0,T) \times \Omega_j$. Accordingly, we define the
  matrices $\mathbf{A}_j$, $\mathbf{B}_{j}$, $\mathbf{M}^{\Gamma_{j}}$,
  $\mathbb{L}_{j,s}^{n,p}$ associated with the finite element method,
  the restriction matrix $\mathbf{Q}_{j,l}$, $\mathbf{Q}_{j,r}$ and
  the solution vector $\mathbf{v}_{j,n}^k$. The subscript "$j$"
    emphasizes the definition of the matrices associated to the
    subdomains  $(0,T)\times \Omega_j$.
\end{remark}

\subsection{Properties of $\mathcal{L}_h$\label{prop_Lh}}


Let us define the fluxes
\begin{displaymath}
  l_{j}^{k}(t,y) = \partial_{\mathbf{n}_j} v^{k}_{j} + S_j v^{k}_{j}, \ r_{j}^{k}(t,y) = \partial_{\mathbf{n}_j} v^{k}_{j} + S_j v^{k}_{j}, \ (t,y) \in (0,T) \times \Omega_j,
\end{displaymath}
for $j=1,2,...,N$ with the two special cases $l_{1}^{k}=r_{N}^{k}=0$
  corresponding to the left and right boundaries and
\begin{displaymath}
  l_{n,j}^{k} = l_{j}^{k}(t_n,y), \ r_{n,j}^{k} = r_{j}^{k}(t_n,y).
\end{displaymath}            
We denote by $\mathbf{l}_{j,n}^k$ (resp. $\mathbf{r}_{j,n}^k$) the nodal interpolation vector of $l_{j,n}^k$ (resp. $r_{j,n}^k$)
\begin{equation}
  \label{rlRobinS2p}
  \quad \left \{ 
    \begin{aligned}
      & \mathbf{r}_{j-1,n}^{k+1} = -\mathbf{l}_{j,n}^k + 2Q_{j,l}
      \cdot \widetilde{S}_j \mathbf{v}_{j,n}^k, \ j=2,3,...,N, \\
      & \mathbf{l}_{j+1,n}^{k+1} = -\mathbf{r}_{j,n}^k + 2
      Q_{j,r} \cdot  \widetilde{S}_j \mathbf{v}_{j,n}^k, \ j=1,2,...,N-1,
    \end{aligned} \right.
\end{equation}
where we have for the Robin transmission condition
\begin{equation*}
  \mathrm{Robin: } \quad  \widetilde{S}_j \mathbf{v}_{j,n}^k = -ip \cdot \mathbf{v}_{j,n}^k,
\end{equation*}
and for the $S_{\mathrm{pade}}^m$ condition
\begin{equation*}
  S_{\mathrm{pade}}^m: \quad  \widetilde{S}_j \mathbf{v}_{j,n}^k  = -i(\sum_{s=0}^m a_s^m) \mathbf{v}_{j,n}^k  + i \sum_{s=1}^m a_s^m d_s^m \bm{\varphi}^{n-1/2}_{j,s}.
\end{equation*}
Then, we define the discrete interface vector $\mathbf{g}$ by 
\begin{equation}
  \label{gdist}
  \mathbf{g} = \lim_{k\rightarrow \infty} \mathbf{g}^k,
\end{equation}
where
\begin{gather*}
  \mathbf{g}^k = (\mathbf{g}_{1}^{k,\top}, \mathbf{g}_{2}^{k,\top}, ..., \mathbf{g}_{N}^{k,\top})^{\top}, \\
  \mathbf{g}_1^k = ( \mathbf{r}_{1,1}^{k,\top}, \mathbf{r}_{1,2}^{k,\top}, \cdots, \mathbf{r}_{1,N_T}^{k,\top})^{\top} \in \mathbb{C}^{N_y \times N_T},\\
  \mathbf{g}_N^k = ( \mathbf{l}_{N,1}^{k,\top}, \mathbf{l}_{N,2}^{k,\top}, \cdots, \mathbf{l}_{N,N_T}^{k,\top})^{\top} \in \mathbb{C}^{N_y \times N_T},\\
  \mathbf{g}_j^k = (\mathbf{l}_{j,1}^{k,\top}, \cdots, \mathbf{l}_{j,N_T}^{k,\top}, \mathbf{r}_{j,1}^{k,\top}, \cdots, \mathbf{r}_{j,N_T}^{k,\top})^{\top}
  \in \mathbb{C}^{2N_y \times N_T},\ j=2,3,...N-1.
\end{gather*}

\begin{proposition}
  \label{InterfGlobalRobin}
  If we consider the Robin transmission condition, the discrete form
  of the interface problem \eqref{InterfacePbContinues} is given
  by \eqref{InterfacePbdisctILgd}. 
\end{proposition}
\proof
  According to \eqref{rlRobinS2p}, and the definitions of $\mathbf{g}^k$, it is easy to verify that
  \begin{equation}
    \label{matXRobin}
    \begin{aligned}
      & X^{j,1} =  -I - 2ip \cdot \mathbf{Q}_{j,l} (\mathbf{A}_j
      - \mathbf{B}_{j})^{-1}
      \mathbf{M}^{\Gamma_{j}} \mathbf{Q}_{j,l}^{\top}, \\
      & X^{j,2} =  - 2ip \cdot \mathbf{Q}_{j,l} (\mathbf{A}_j
      - \mathbf{B}_{j})^{-1}
      \mathbf{M}^{\Gamma_{j}} \mathbf{Q}_{j,r}^{\top}, \\
      & X^{j,3} =  - 2ip \cdot \mathbf{Q}_{j,r} (\mathbf{A}_j
      - \mathbf{B}_{j})^{-1}
      \mathbf{M}^{\Gamma_{j}} \mathbf{Q}_{j,l}^{\top}, \\
      & X^{j,4} =  -I - 2ip \cdot \mathbf{Q}_{j,r} (\mathbf{A}_j
      - \mathbf{B}_{j})^{-1}
      \mathbf{M}^{\Gamma_{j}} \mathbf{Q}_{j,r}^{\top},
    \end{aligned}
  \end{equation}  
  and
  \begin{equation}
    \label{dRobin}
    \begin{split}
      \mathbf{d}_{j,l} & = 2ip \cdot \mathbf{Q}_{j-1,r} (\mathbf{A}_{j-1} - \mathbf{B}_{j-1})^{-1} \mathbf{F}_{j-1}, \ j=2,3,...,N, \\
      \mathbf{d}_{j,r} & = 2ip \cdot \mathbf{Q}_{j+1,l} (\mathbf{A}_{j+1} - \mathbf{B}_{j+1})^{-1} \mathbf{F}_{j+1}, \ j=1,2,...,N-1.
    \end{split}
  \end{equation}
  

\begin{proposition}
  \label{InterfGlobalS2P}
  If we consider the transmission condition $S_{\mathrm{pade}}^m$,
  then \eqref{InterfacePbdisctILgd} is the discrete form of the
  interface problem \eqref{InterfacePbContinues}.
\end{proposition}
\proof
  By using \eqref{rlRobinS2p}, we have
  \begin{displaymath}
    \widetilde{S}_j \mathbf{v}_{j,n}^k  = - c_a \mathbf{v}_{j,n}^k  + i \sum_{s=1}^m a_s^m d_s^m \sum_{p=1}^{n} \mathbb{L}_{j,s}^{n,p} \mathbf{v}_{j,p}^k = - c_a \mathbf{v}_{j,n}^k  + i\sum_{p=1}^{n} \big( \sum_{s=1}^m a_s^m d_s^m \mathbb{L}_{j,s}^{n,p} \big) \mathbf{v}_{j,p}^k.
  \end{displaymath}
  We could easily verify that  
  \begin{equation}
    \label{matXSP}
    \begin{aligned}
      & X^{j,1} =  -I + 2 \mathbf{Q}_{j,l} \mathbf{B}_j^S  (\mathbf{A}_j
      - \mathbf{B}_{j})^{-1}
      \mathbf{M}^{\Gamma_{j}} \mathbf{Q}_{j,l}^{\top}, \\
      & X^{j,2} =  2\mathbf{Q}_{j,l} \mathbf{B}_j^S  (\mathbf{A}_j
      - \mathbf{B}_{j})^{-1}
      \mathbf{M}^{\Gamma_{j}} \mathbf{Q}_{j,r}^{\top}, \\
      & X^{j,3} =   2\mathbf{Q}_{j,r} \mathbf{B}_j^S  (\mathbf{A}_j
      - \mathbf{B}_{j})^{-1}
      \mathbf{M}^{\Gamma_{j}} \mathbf{Q}_{j,l}^{\top}, \\
      & X^{j,4} = -I + 2\mathbf{Q}_{j,r} \mathbf{B}_j^S (\mathbf{A}_j
      - \mathbf{B}_{j})^{-1}
      \mathbf{M}^{\Gamma_{j}} \mathbf{Q}_{j,r}^{\top},
    \end{aligned} 
  \end{equation}
  and
  \begin{equation}
    \label{dSP}
    \begin{split}
      \mathbf{d}_{j,l} & = \mathbf{Q}_{j-1,r} \mathbf{B}_{j-1}^S (\mathbf{A}_{j-1} - \mathbf{B}_{j-1})^{-1} \mathbf{F}_{j-1}, \ j=2,3,...,N, \\
      \mathbf{d}_{j,r} & = \mathbf{Q}_{j+1,l} \mathbf{B}_{j+1}^S  (\mathbf{A}_{j+1} - \mathbf{B}_{j+1})^{-1} \mathbf{F}_{j+1}, \ j=1,2,...,N-1,
    \end{split}
  \end{equation}
  where the matrix $\mathbf{B}_j^S$ is defined by
  \begin{displaymath}
    \mathbf{B}_j^S = 
    \begin{pmatrix}                                  
      -c_a I +  \displaystyle i\sum_{s=1}^m a_s^m d_s^m \mathbb{L}_{j,s}^{1,1} \\
      \displaystyle i\sum_{s=1}^m a_s^m d_s^m \mathbb{L}_{j,s}^{2,1} & \displaystyle  -c_a I + i\sum_{s=1}^m a_s^m d_s^m \mathbb{L}_{j,s}^{2,2} \\
      \vdots & \vdots  & \ddots \\
      \displaystyle i\sum_{s=1}^m a_s^m d_s^m \mathbb{L}_{j,s}^{N_T,1} & \displaystyle i\sum_{s=1}^m a_s^m d_s^m \mathbb{L}_{j,s}^{N_T,2} & \cdots & -c_a I + \displaystyle i\sum_{s=1}^m a_s^m d_s^m \mathbb{L}_{j,s}^{N_T,N_T}
    \end{pmatrix}.
  \end{displaymath}
  %

Let us now study the structure of the subblock of $\mathcal{L}_h$
  for time independent potential $ V = V(x,y)$. More specifically, we
  focus on
\begin{equation}
  \label{blockinf}
  \begin{aligned}
    & X^{1,4} = \{ x^{1,4}_{n,s} \}_{1 \leqslant n,s \leqslant N_T}, \\
    & X^{j,1} = \{ x^{j,1}_{n,s} \}_{1 \leqslant n,s \leqslant N_T}, X^{j,2}=\{ x^{j,2}_{n,s} \}_{1 \leqslant n,s \leqslant N_T},\\
    & X^{j,3} = \{ x^{j,3}_{n,s} \}_{1 \leqslant n,s \leqslant N_T}, X^{j,4}=\{ x^{j,4}_{n,s} \}_{1 \leqslant n,s \leqslant N_T},j=2,3,...,N-1,\\
    & X^{N,1} = \{ x^{N,1}_{n,s} \}_{1 \leqslant n,s \leqslant N_T}.
  \end{aligned}
\end{equation}
where $x^{j,1}_{n,s},x^{j,2}_{n,s},x^{j,3}_{n,s},x^{j,4}_{n,s} \in \mathbb{C}^{N_y \times N_y}$ are submatrices.

\begin{figure}[!htbp]
  \centering
  \includegraphics[width=0.3\textwidth]{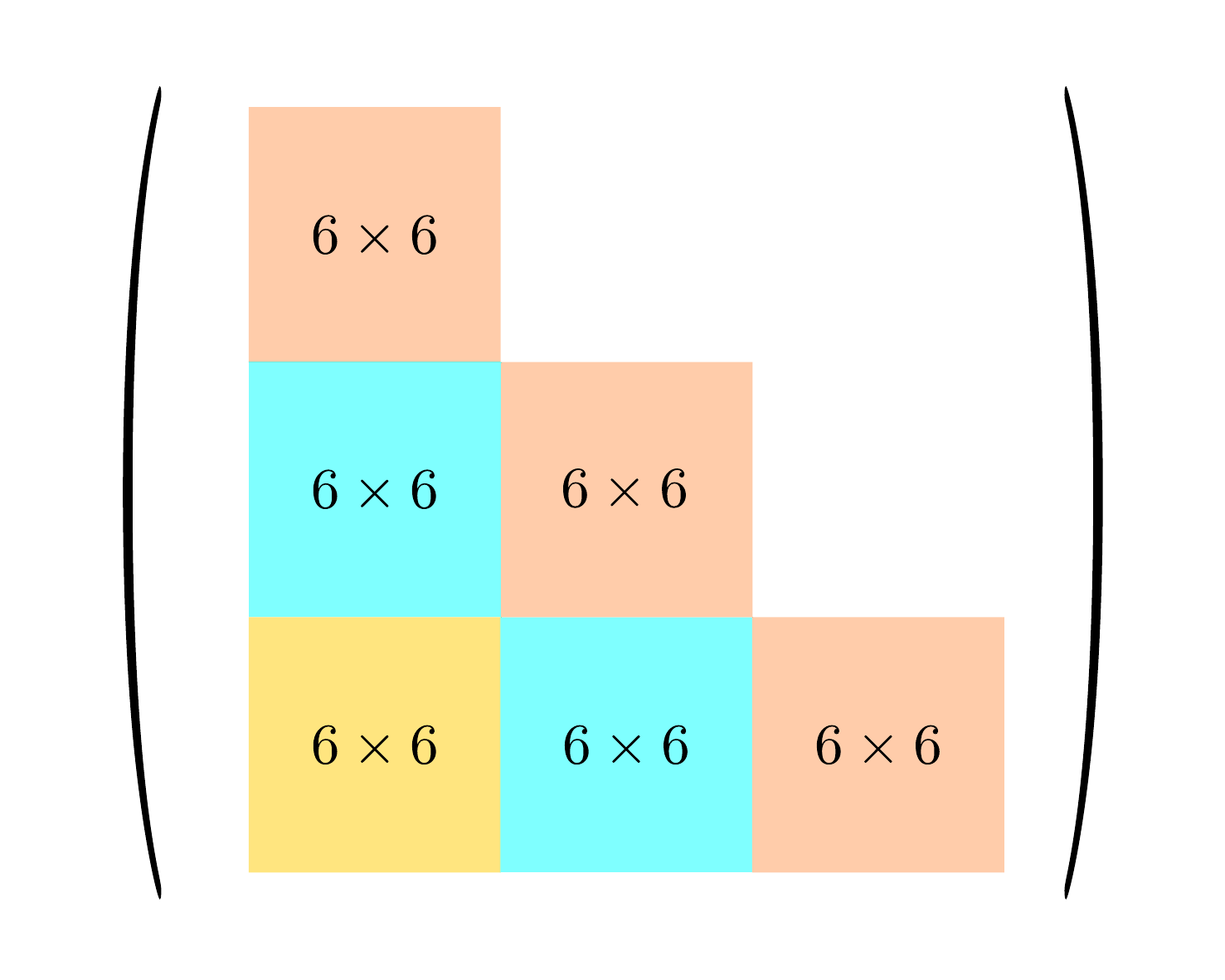}
  \caption{Block structure, $N_T=3$, $N_y=6$.}
  \label{matblock}
\end{figure}

Each subblock $X_j^{1,2,3,4}$ are made of submatrices which are set on
  very specific positions. This structure is presented in Figure
  \ref{matblock} for 3 time steps and 6 nodes on the interface between
  two subdomains. We see that each sub-diagonal block is identical. We present
this property mathematically in proposition \ref{Prop_Blockprop}
with the two transmission condition Robin or $S_{\mathrm{pade}}^m$. The
demonstration is similar to the one obtained for one dimensional Schr\"{o}dinger
equation \cite{XF20151d}. The formal difference between
  dimension one and dimension two is that the flux are
  scalar in one dimension and vectors in two dimensions.
%

\begin{proposition}
  \label{Prop_Blockprop}
  For the transmission Robin (resp. $S_{\mathrm{pade}}^m$), assuming
  that the linear system \eqref{NproblemRobin}
  (resp. \eqref{NproblemSP}) is not singular, if the potential is
    time independent $V=V(x,y)$, then the
  matrices $X^{1,4}$, $X^{j,1}$,$X^{j,2}$, $X^{j,3}$, $X^{j,4}$,
  $j=2,3,...,N-1$ and $X^{N,1}$ are block lower triangular matrices
  and they satisfy 
  \begin{displaymath}
    \begin{split}
      & x^{1,4}_{n,s} = x^{1,4}_{n-1,s-1},\\
      & x^{j,1}_{n,s} = x^{j,1}_{n-1,s-1},\
      x^{j,2}_{n,s} = x^{j,2}_{n-1,s-1},\\
      & x^{j,3}_{n,s} = x^{j,3}_{n-1,s-1},\
      x^{j,4}_{n,s} = x^{j,4}_{n-1,s-1}, j=2,3,...,N-1,\\
      & x^{N,1}_{n,s} = x^{N,1}_{n-1,s-1},
    \end{split}
  \end{displaymath} 
  for $2 \leqslant s \leqslant n \leqslant N_T $.
\end{proposition}

If the potential is constant, the sub-blocks of $\mathcal{L}_h$ are
  simpler and follow the following proposition.

\begin{proposition}
  \label{Prop_V0}
  Assuming that the matrix $\mathbf{A}_j  - \mathbf{B}_{j}$ is not
  singular with the Robin or the $S_{\mathrm{pade}}^m$ transmission
  condition, and that the mesh is uniform and
  the size of subdomains $\Omega_j$ are equal, if the potential $V$ is
  constant, then the subblocks of
  $\mathcal{L}_h$ satisfy 
  \begin{equation}
    \label{blockidt}
    \begin{aligned}
      & X^{2,1} = X^{3,1} = \cdots = X^{N,1}, \quad \ X^{2,2} = X^{3,2} = \cdots = X^{N-1,2},\\
      & X^{2,3} = X^{3,3} = \cdots = X^{N-1,3}, \ X^{1,4} = X^{2,4} = \cdots =   X^{N-1,4}.
    \end{aligned}
  \end{equation}
\end{proposition}

\proof
  Thanks to the hypothesis of the proposition, the geometry of each
    subdomain is identical. Thus, the various matrices coming from the
    assembly of the finite element methods are the same. Therefore, we have 
  \begin{displaymath}
    \mathbb{A}_1 = \mathbb{A}_2=\cdots=\mathbb{A}_N, \
    \mathbb{M}^{\Gamma_1}=\mathbb{M}^{\Gamma_2} = \cdots
    =\mathbb{M}^{\Gamma_N}, 
  \end{displaymath}
  and the restrictions matrices satisfy
  \begin{displaymath}
    Q_{1,l} = Q_{2,l} = \cdots = Q_{N,l}, \ Q_{1,r} = Q_{2,r} = \cdots = Q_{N,r}.
  \end{displaymath}
  Since $V$ is constant, then $\mathbb{M}_{1,W_{n}}=\mathbb{M}_{2,W_{n}}=...=\mathbb{M}_{N,W_{n}}$. Thus by definitions, we have
  \begin{displaymath}
    \mathbf{A}_1 = \mathbf{A}_2=\cdots=\mathbf{A}_N, \
    \mathbf{B}_{1}=\mathbf{B}_{2} = \cdots
    =\mathbf{B}_{N},\
    \mathbf{Q}_{1} = \mathbf{Q}_{2} = \cdots = \mathbf{Q}_{N},
  \end{displaymath}
  and $\mathbf{B}_{1}^S=\mathbf{B}_{2}^S = \cdots =\mathbf{B}_{N}^S$ for
  the transmission condition $S_{\mathrm{pade}}^m$. The conclusion
  therefore follows from \eqref{matXRobin} and  \eqref{matXSP}. 

\section{New algorithms for the
  Schr\"{o}dinger equation with a constant
  potential \label{new_algo_sec}}

Thanks to the analysis yielded in previous section, we can build explicitly the
  interface problem \eqref{InterfacePbdisctILgd}. 
  The main idea of our new algorithm is therefore to explicitly construct the matrix
  $\mathcal{L}_h$ and the vector $d$. Following propositions
  \eqref{Prop_Blockprop} and \eqref{Prop_V0}, it is sufficient to compute four subblocks to
explicitly build the matrix $\mathcal{L}_h$. Without loss of
generality, we construct the blocks $X^{2,1}$, $X^{2,2}$, $X^{2,3}$
and $X^{2,4}$. Furthermore, only the first $N_y$ columns of each block
are necessarily computed.

\subsection{New SWR algorithm}

We propose here a new SWR algorithm for the Schr\"{o}dinger equation
with a constant potential $V$. We always suppose that the size of each
subdomain is identical and the mesh is uniform. 
\begin{center}
  \begin{minipage}{0.85\textwidth}
    \begin{algorithm}[H]
      \caption{New SWR algorithm, constant potential $V$ }
      \label{Algo_newSWR}
      \begin{algorithmic}[1] 
        \STATE Build $\mathcal{L}_h$ and $d$ in relation $(I-\mathcal{L}_h) \mathbf{g} = d$ explicitly.
        \STATE Solve the linear system $(I-\mathcal{L}_h) \mathbf{g} =d$ by an iterative method (fixed point method or Krylov methods).
        \STATE Solve the Schr\"{o}dinger equation on each subdomain on $(0,T) \times \Omega_j$ using the flux obtained from step 2.
      \end{algorithmic}
    \end{algorithm}
  \end{minipage}
\end{center}
Mathematically, this new SWR algorithm is identical to the classical
one. The main novelty here is that we construct
explicitly the matrix $\mathcal{L}_h$ and the vector $d$
(corresponding to the first
step), while in the classical algorithm, $\mathcal{L}_h$ remains an abstract operator. It is not usual to build explicitly such
  huge operator, but as we will
  see, its computation is not costly.

We show below the construction of $\mathcal{L}_h$ and $d$ for the Robin transmission
condition. This is is based on the formulas \eqref{matXRobin} and
\eqref{dRobin}. For the transmission condition $S_{\mathrm{pade}}^m$,
the idea is similar, but involves \eqref{matXSP} and \eqref{dSP}.  According to the proposition \ref{InterfGlobalRobin}, the
column $s$ of $X^{2,1}$ and $X^{2,3}$ are  
\begin{displaymath}
  \begin{split}
    X^{2,1} \mathbf{e}_s & = -\mathbf{e}_s - 2ip \cdot \mathbf{Q}_{2,l} \mathbf{M}^{\Gamma_{2}} \mathbf{Q}_{2,l}^{\top} \mathbf{M}^{\Gamma_{2}} \mathbf{Q}_{2,l}^{\top} \mathbf{e}_s, \\
    X^{2,3} \mathbf{e}_s & = - 2ip \cdot \mathbf{Q}_{2,r} (\mathbf{A}_2
    - \mathbf{B}_{2})^{-1} \mathbf{M}^{\Gamma_{2}} \mathbf{Q}_{2,l}^{\top} \mathbf{e}_s,
  \end{split}
\end{displaymath}                               
where the vector $\mathbf{e}_s =(0,0,...,1,...0) \in \mathbb{C}^{N_T
  \times N_y}$, all its elements being zero to the exception of the $s$-th, which is
one. The element $\mathbf{M}^{\Gamma_{2}} \mathbf{Q}_{2,l}^{\top} \mathbf{e}_s$ being
  a vector, it is necessary to compute one time the application of
  $(\mathbf{A}_2 - \mathbf{B}_{2})^{-1}$ to a vector. Similarly, we have 
\begin{displaymath}
  \begin{split}
    X^{2,2} \mathbf{e}_s & = - 2ip \cdot \mathbf{Q}_{2,l} (\mathbf{A}_2
    - \mathbf{B}_{2})^{-1} \mathbf{M}^{\Gamma_{2}} \mathbf{Q}_{2,r}^{\top} \mathbf{e}_s, \\
    X^{2,4} \mathbf{e}_s & = -\mathbf{e}_s - 2ip \cdot \mathbf{Q}_{2,r} (\mathbf{A}_2
    - \mathbf{B}_{2})^{-1}
    \mathbf{M}^{\Gamma_{2}} \mathbf{Q}_{2,r}^{\top} \mathbf{e}_s,
  \end{split}
\end{displaymath} 
In conclusion, to know the first $N_y$ columns of $X^{2,1}$,
$X^{2,2}$, $X^{2,3}$ and $X^{2,4}$, it is sufficient to compute $2N_y$
times the application of $(\mathbf{A}_2 - \mathbf{B}_{2})^{-1}$ to a
vector. In other words, this amounts to solve the Schr\"{o}dinger equation on a single
subdomain $2N_y$ times to build the
matrix $\mathcal{L}_h$. {Without loss of generality, since the
  geometry of each sub domain is identical and the matrix
  $\mathcal{L}_h$ is independent of the initial solution, we focus
  here on the domain $(0,T) \times \Omega_2$.} The resolutions being all independent, we
  can solve them on different processors using MPI paradigm. We fix
one MPI process per domain. To construct the matrix $\mathcal{L}_h$,
we use the $N$ MPI processes to solve the equation on a single
subdomain ($(0,T) \times \Omega_2$) $2N_y$ times. Each MPI process
therefore solves the Schr\"{o}dinger equation on a single
subdomain maximum  
%

%
\begin{displaymath}
  N_{\mathrm{mpi}} := [\frac{2N_y}{N}] + 1 \ \text{times,}
\end{displaymath}
where $[x]$ denotes the integer part of $x$. This construction is
therefore super-scalable in theory. Indeed, if $N$ is doubled, then the size of
the subdomains is divided by two and $N_{\mathrm {mpi}} $ is also
approximately halved.  

Concerning the computational phase, the transposed matrix of
$\mathcal{L}_h$ is stored in a distributed manner using the
PETSc library \cite{petsc-user-ref}. As shown by \eqref{matL}, the
first block column of $\mathcal{L}_h$ lies in MPI process 0. The second
and third blocks columns are in MPI process 1, and so on for other
processes. The size of each block is $(N_T \times N_y) \times (N_T
\times N_y)$. Each block contains $(N_T+1) \times N_T/2 \times N_y^2$
nonzero elements.                           

The construction of the vector $\mathbf{d}$ is similar. According to
\eqref{dRobin}, one needs to apply $(\mathbf{A}_j -
\mathbf{B}_{j})^{-1}$ to vector $\mathbf{F}_j$ for $j=1,2,...,N$. In
other words, we solve the Schr\"odinger equation on each subdomain
one time. Again, the vector is stored in a distributed manner using the PETSc
library. 

\begin{minipage}{10cm}
  \centering
  \begin{displaymath}
    \mathbf{d} = 2 \begin{pmatrix}
      \mathbf{d}_{1,r} \\
      \vdots \\
      \mathbf{d}_{j,l} \\
      \mathbf{d}_{j,r} \\
      \vdots \\
      \mathbf{d}_{N,l}
    \end{pmatrix}
    \begin{array}{c} 
      \rdelim \}{1}{4pt}[MPI 0, $N_T\times N_y$ elements] \\
      \\
      
      \rdelim \}{2}{4pt}[MPI $j$, $2N_T\times N_y$ elements] \\
      \\
      \\
      \rdelim \}{1}{4pt}[MPI $N$, $N_T\times N_y$ elements] \\
    \end{array}
  \end{displaymath} 
\end{minipage}

\subsection{New DDS algorithm}

Since we have interpreted the DDS method as sequence of some SWR
methods in Algorithm \ref{AlgoDDS_swr}, we can apply directly the ideas developped in the previous sections to the DDS method. Let us denote the interface problem of 
\begin{displaymath}
  \left \{ \begin{array}{ll}
      i\partial_t u + \Delta u + V u = 0, &(t,x,y) \in (t_{n-1},t_n) \times \Omega,\\
      u(0,x,y) = u_{n-1}(x,y), &(x,y) \in \Omega, \\
    \end{array} \right. 
\end{displaymath}	
by 
\begin{equation}
\label{InterfaceDDS}
(I-\mathcal{L}_{h,n})\mathbf{g}_{n}=d_{n},
\end{equation}
where $\mathcal{L}_{h,n}$ is interface matrix and $d_n$ is vector.  In
the classical algorithm, $\mathcal{L}_{h,n}$ remains an abstract
operator. We propose to build it explicitly in the new DDS
algorithm. Actually, thanks to the following proposition, the
computation of the complete operator $\mathcal{L}_{h,n}$ only require
to compute $\mathcal{L}_{h,1}$. 
\begin{proposition}
  For the transmission Robin and $S_{\mathrm{pade}}^m$, the interface matrix satisfies
  \begin{displaymath}
    \mathcal{L}_{h,1} = \mathcal{L}_{h,2}= ... = \mathcal{L}_{h,N_T}.
  \end{displaymath}
\end{proposition}
\proof
According to \eqref{matXRobin} and \eqref{matXSP}, the interface matrix is independent of the initial datum, thus the conclusion follows directly.
%
\begin{center}
  \begin{minipage}{0.85\textwidth}
    \begin{algorithm}[H]
      \caption{New DDS algorithm, constant potential $V$}
      \label{Algo_newDDS}
      {\footnotesize 1:} Build $\mathcal{L}_{h,1}$ explicitly.\\
      {\footnotesize 2:} The initial datum is $u_{0}$. \\
      \For{$n=1,2,...,N_T$}{
        {\footnotesize 2.1:} Build $d_{n}$ on time $t_n$,\\
        {\footnotesize 2.2:} Solve the linear system $(I-\mathcal{L}_{h,n})\mathbf{g}_{h,n}=d_{n}$ by an iterative method, where the initial vector is chosen as $\mathbf{g}_{h,n-1}$.\\
        {\footnotesize 2.3:} Solve the Schr\"{o}dinger equation on each subdomain $(t_{n-1},t_n) \times \Omega_j$ using the flux from step 2.2 to compute $u_{n}$. 
      }
    \end{algorithm}
  \end{minipage}
\end{center}
Mathematically, this new DDS algorithm is identical to the classical one. Compared with the new SWR algorithm, the construction of the interface matrix is less costly since the size of $\mathcal{L}_{h,1}$ is smaller than that of $\mathcal{L}_{h}$.
%

\section{Preconditioned algorithms for general linear potential}

In the case of a non constant 
potential, the proposition \ref{Prop_V0} does not hold. Thus we could
not construct easily the matrix $\mathcal{L}_h$. The aim of this
section is to present preconditioned algorithms for the
Schr\"{o}dinger equation \eqref{Schrodingereq} with a non constant 
potential $V=V(t,x,y)$. Adding a preconditioner to the equation
\eqref{InterfacePbdisctILgd} leads to a preconditioned SWR algorithm 
\begin{equation}
  P^{-1} (I - \mathcal{L}_h) \mathbf{g} = P^{-1} d,
\end{equation}
We propose here to use
\begin{displaymath}
  P = I - \mathcal{L}_0,
\end{displaymath}
where $\mathcal{L}_0$ denotes the interface matrix in
\eqref{InterfacePbdisctILgd} defined for the free Schr\"{o}dinger
equation when $V=0$. As mentioned in the previous section, it could be
easily constructed numerically and it is stored in a distributed
manner. 
%
For any vector $y$, the vector $x: = P^{-1} y$ is computed by solving the linear system
\begin{displaymath}
  P x = y.
\end{displaymath}  
with Krylov methods (for example GMRES or BiCGStab). However, the size
of $P$ increases linearly with the number of subdomains $N$. Also the
number of involved operations for multiplying $\mathcal {L}_0$ and a vector
(which is the basic operation of Krylov method) is not negligible
compared to the computation of the solutions to the equations on each subdomain if $N$ is
large. Thus, the application of the preconditioner increases. 

We could derive straightforwardly a preconditioned DDS algorithm from
the point of view of Algorithm \ref{AlgoDDS_swr}. The preconditioner
for all time steps $n=1,2,...,N_T$ is always chosen as  
\begin{displaymath}
  P=I - \mathcal{L}_{0,1},
\end{displaymath}
where $\mathcal{L}_{0,1}$ is the interface matrix of \eqref{InterfaceDDS} defined for $V=0$ and $n=1$.
%

\begin{remark}
  The idea of these new algorithm are not limited to the Schr\"odinger
  equation. It could also be applied to some other PDEs, like heat
  equation etc. The only limitation is that the mesh and the decomposition should be uniform.
\end{remark}

\section{Numerical results}
\label{Sec_NumResults}

\begin{figure}[!htbp]
  \centering
  \includegraphics[resolution=600]{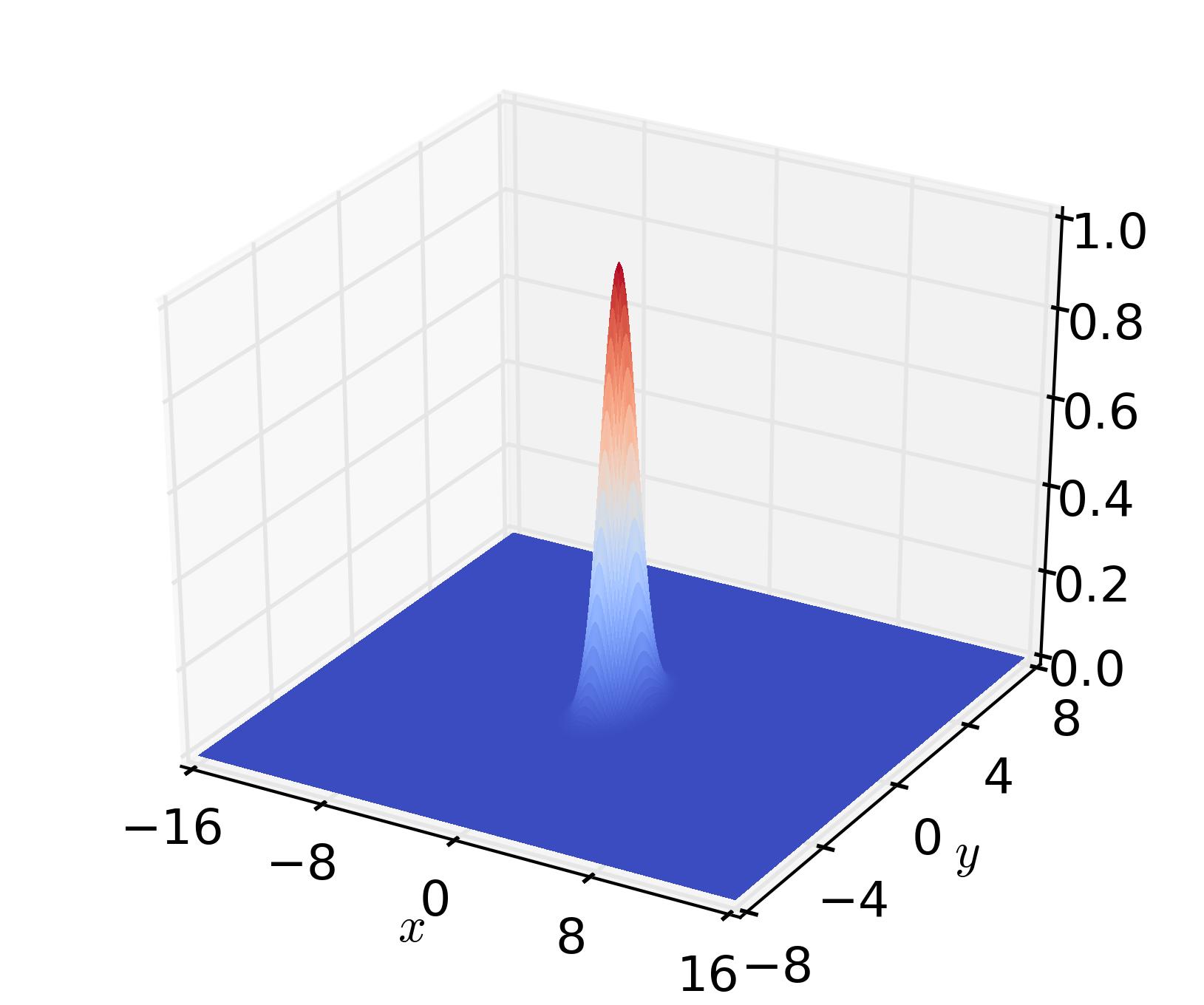} 
  \caption{Initial datum $|u_0|$, $u_0 = e^{-x^2-y^2 - 0.5i x}$.}
  \label{Solinit0}
\end{figure}

The complete domain $\Omega=(-16,16) \times (-8,8)$ is decomposed into
$N$ equal subdomains. The
size of a single cell is $\Delta x \times \Delta y$. We consider two
different meshes  
\begin{gather*}
  \Delta x=1/128, \ \Delta y=1/8;\\ \Delta x=1/2048, \ \Delta y=1/128.
\end{gather*}
With the first mesh, it is possible to solve the Schr\"{o}dinger
equation \eqref{Schrodingereq} on the entire domain on a single node
of a cluster composed of 92 nodes (16 cores/node, Intel Sandy Bridge
E5-2670, 32GB memory/node). Thus we could observe if the parallel
algorithms allow to reduce the total computation time of the
sequential algorithm. We are interested in the strong scalability up to 1024 subdomains. 
The initial datum in this section is 
\begin{equation}
  \label{Initdata}
  u_0(x,y) = e^{-x^2-y^2 - 0.5i x}. \quad \text{(see Figure \ref{Solinit0})}
\end{equation}

This section is composed of two subsections. The first one is devoted to the free
Schr\"{o}dinger equation ($V=0$). In the second, we consider the
Schr\"{o}dinger equation  with a linear potential $V=x^2+y^2$. 

\begin{remark}
  \label{rmk1step}
For the DDS algorithm, we consider mostly the convergence properties
  (for example the number of iterations) of the first time step. In other words,
  we only study the evolution between $t^0 \rightarrow t^1$ to compute $v^1$ with the
  optimized Schwarz method.
\end{remark}

\begin{remark}
  \label{rmkrobin}
  The theoretical optimal parameter $p$ in the Robin transmission
  condition is not at hand for us. We then seek for the best parameter
  numerically for each case.  
\end{remark}

\subsection{The free Schr\"{o}dinger equation}
\label{Sec_FreeSch}
In this part, we first compare the SWR method and the DDS
algorithm. The DDS algorithm being more efficient, we next compare
  the classical DDS and the new DDS algorithms. Finally, the
influence of parameters is studied in 
Section \ref{Sec_InfluenceParameter}. 
\subsubsection{SWR vs. DDS}

We compare the SWR and the DDS methods in the framework of the
classical algorithm. The time step is set to $\Delta t = 0.01$. The
mesh size is $\Delta x = 1/128$, $\Delta y = 1/8$. Both methods are applied
on time-space domain $(0,T)\times \Omega$ with various final time
$T$. We denote by $N_{\mathrm{iter}}^{S}$ the number of iterations
required to obtain convergence of the SWR method, $N_{\mathrm{iter}}^{D}$ the number of
iterations of the first time step of the DDS method and $T^S$
(resp. $T^D$) the total computation time of the SWR method (resp. DDS
method). Table \ref{SWRDDSCompare} and Table \ref{SWRDDSCompareRobin}
present the numbers of iterations and the computation times of both
methods for $N =2$ and $N = 32$, where the transmission conditions are
$S_{\mathrm{pade}}^m, m =5$ and Robin respectively. The initial vector
is the zero vector and the GMRES method is used on the interface
problem.  We could see that  $N_{\mathrm{iter}}^{S}$ increases with the
final time $T$ and $N_{\mathrm{iter}}^{S}>N_{\mathrm{iter}}^{D}$. Thus
$T^S > T^D$. 



%
\begin{table}[!htbp] 
  \centering
  \footnotesize
  \renewcommand{\arraystretch}{1.2}
  \caption{Number of iterations and computation time of the SWR method and the DDS method for $N=2,32$ with the transmission condition $S_{\mathrm{pade}}^5$.}
  \begin{tabular}{|c|c|c|c|c|c|c|c|c|}
    \hline
    $T$ & \multicolumn{4}{c|}{$N=2$} & \multicolumn{4}{c|}{$N=32$}  \\
    \hline
    & $N_{\mathrm{iter}}^{S}$ & $N_{\mathrm{iter}}^{D}$ & $T^S$ & $T^D$ & $N_{\mathrm{iter}}^{S}$ & $N_{\mathrm{iter}}^{D}$ & $T^S$ & $T^D$ \\
    \hline 
    0.05 & 17 & 9 & 17.6 & 12.3 & 17 & 10 & 1.5 & 1.9 \\
    0.1 & 25 & 9 & 44.5 & 21.5 & 25 & 10 & 3.5 & 1.8 \\
    0.15 & 30 & 9 & 78.4 & 30.9 & 31 & 10 & 6.4 & 2.6 \\
    0.2 & 44 & 9 & 147.9 & 40.0 & 45 & 10 & 11.8 & 3.4 \\
    0.25 & 51 & 9 & 215.0 & 49.3 & 52 & 10 & 16.9 & 4.1 \\
    0.3 & 55 & 9 & 271.0 & 58.6 & 55 & 10 & 21.2 & 4.8 \\
    0.35 & 58 & 9 & 332.1 & 67.9 & 59 & 10 & 26.3 & 5.6 \\
    0.4 & 61 & 9 & 402.9 & 77.2 & 62 & 10 & 32.0 & 6.3 \\
    0.45 & 64 & 9 & 474.5 & 86.4 & 65 & 10 & 37.6 & 7.1 \\
    0.5 & 68 & 9 & 557.6 & 96.1 & 73 & 10 & 46.3 & 7.8\\
    \hline
  \end{tabular}
  \label{SWRDDSCompare}
\end{table}
\begin{table}[!htbp] 
  \centering
  \footnotesize
  \renewcommand{\arraystretch}{1.2}
  \caption{Number of iterations and computation time of the SWR method and the DDS method for $N=2,32$ with Robin transmission condition.}
  \begin{tabular}{|c|c|c|c|c|c|c|c|c|}
    \hline
    $T$ & \multicolumn{4}{c|}{$N=2$} & \multicolumn{4}{c|}{$N=32$}  \\
    \hline
    & $N_{\mathrm{iter}}^{S}$ & $N_{\mathrm{iter}}^{D}$ & $T^S$ & $T^D$ & $N_{\mathrm{iter}}^{S}$ & $N_{\mathrm{iter}}^{D}$ & $T^S$ & $T^D$ \\
    \hline 
    0.05 & 19 & 11 & 13.0 & 15.9 & 20 & 11 & 1.1 & 0.6 \\
    0.1 & 28 & 11 & 32.6 & 16.8 & 29 & 11 & 2.7 & 0.6 \\
    0.15 & 44 & 11 & 73.3 & 17.9 & 47 & 11 & 6.1 & 0.7 \\
    0.2 &  57 & 11 & 123.1 & 18.9 & 57 & 11 & 9.6 & 0.8 \\
    0.25 & 76 & 11 & 204.6 & 19.9  & 80 & 11 & 16.6 & 0.9 \\
    0.3 &  85 & 11 & 271.4 & 21.0 & 89 & 11 & 22.1 & 1.0 \\
    0.35 & 90 & 11 & 333.6 & 21.9 & 92 & 11 & 26.7 & 1.0 \\
    0.4 &  100 & 11 & 425.3 & 23.0 & 102 & 11 & 33.5 & 1.2 \\
    0.45 & 109 & 11 & 519.4 & 24.0 & 208 & 11 & 40.0 & 1.3 \\
    0.5 & 114 & 11 & 601.2 & 25.0 & 116 & 11 & 47.3 & 1.3 \\
    \hline
  \end{tabular}
  \label{SWRDDSCompareRobin}
\end{table}
In the next subsections, without special statement, we only consider the DDS method since it does requires less computation time.

\subsubsection{Comparison of classical and new algorithms}

In this part, we are interested in the performance (number of
iterations and computation time) of the classical and the new
algorithms with the two transmission conditions. We observe the strong
scalability of the two algorithms. Both algorithms and transmission
conditions are compared in the framework of the DDS method. The final
time is $T = 0.5$ and the time step is fixed as $\Delta t = 0.01$. The
GMRES method is used on the interface problem. The initial vector is
the zero vector. Since there is no theoretical result for us on the
choice of the optimal parameter $p$ in the Robin transmission
condition, we make tests with different $p$ to find the numerically
optimal one. However, we will see in the next subsection that the
number of iterations and the computing time are not sensitive to $p$
using the GMRES method on the interface problem. Thus, it is difficult
to choose the optimal one. We take $p = 15$ for the mesh $ \Delta x =
1/128$, $\Delta y = 1/8$ and $ p = 10$ for the mesh $\Delta x =
1/2048$, $\Delta y = 1/64$. For the transmission condition
$S_{\mathrm{pade}}^m$, we take $m = 5$, the numerical optimal according to some tests. 

%

\begin{figure}[!htbp]
  \centering
  \includegraphics[width=0.46\textwidth]{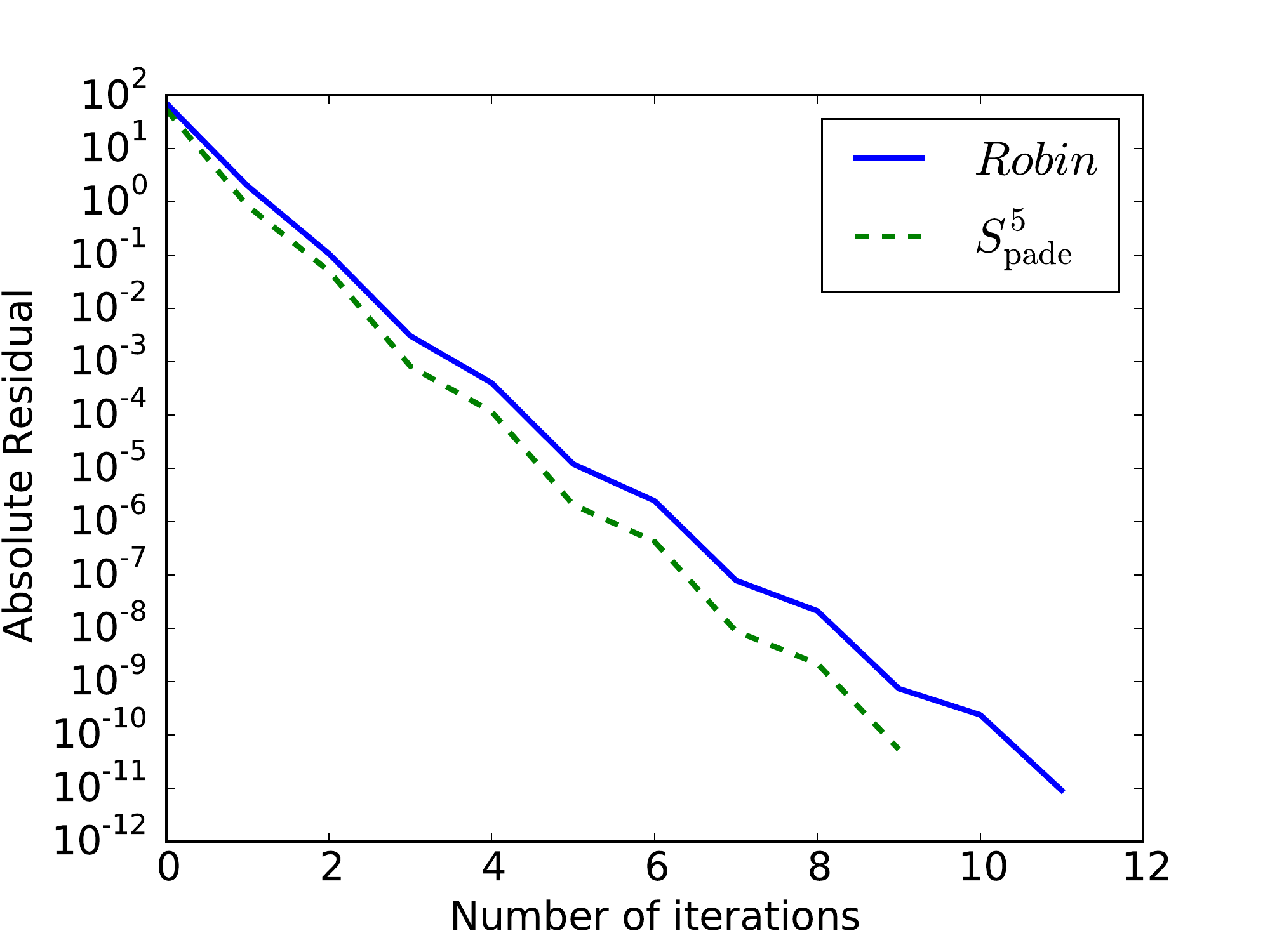}
  \includegraphics[width=0.46\textwidth]{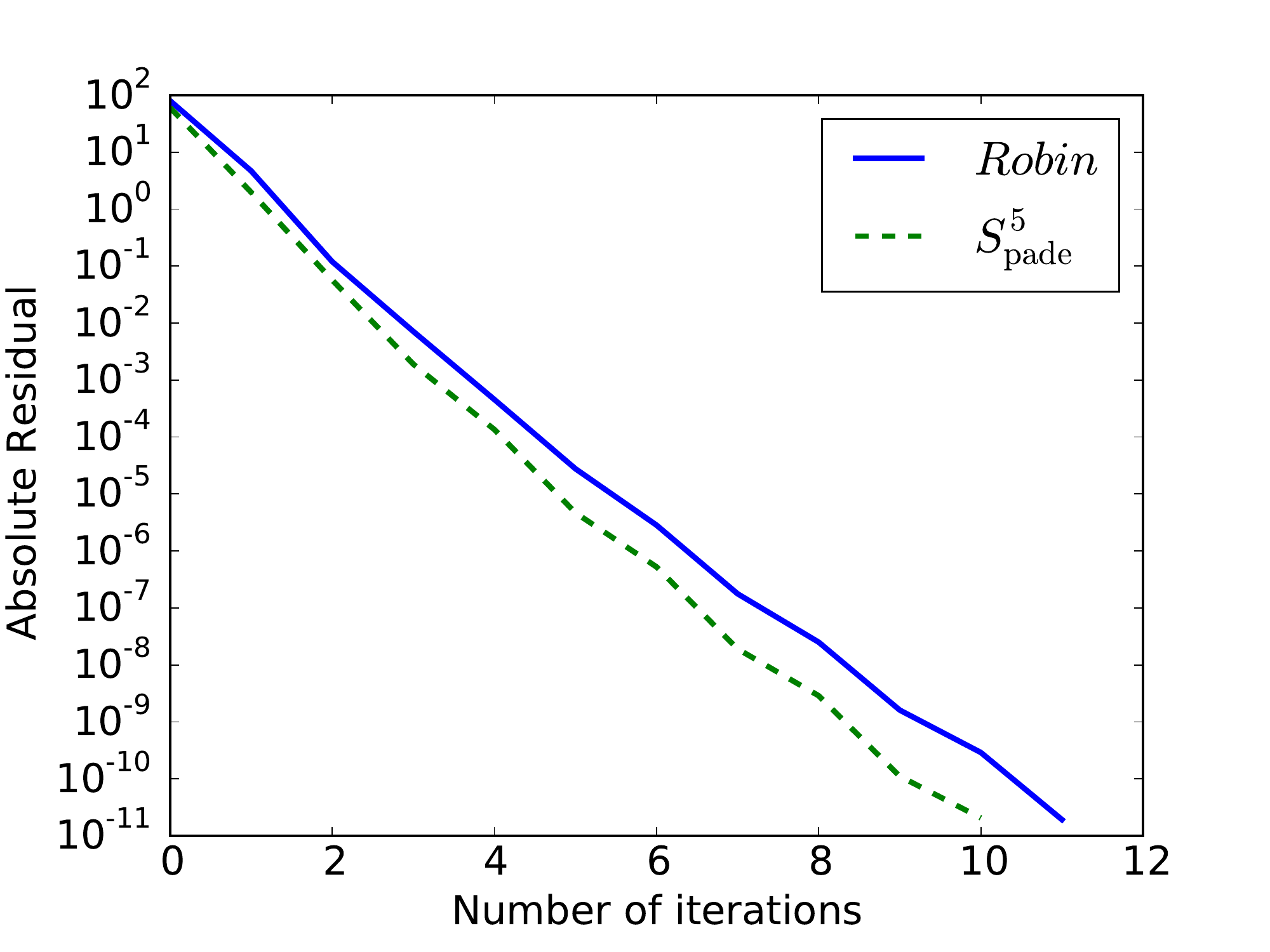}
  \caption{Convergence history of the first time step for $N=2$ (left) ,$N=32$ (right).}
  \label{Hist_2_32_V0}
\end{figure}

We first consider the mesh $\Delta x = 1/128$, $\Delta y =
1/8$. Figure \ref{Hist_2_32_V0} presents the convergence history of
the first time step for $N = 2$ and $N = 32$. We also show the number
of iterations of the first time step and the total computation time in
Table \ref{Time_clnew_T>dt_128_V0} for $N = 2,4,8,16,32$. The boundary conditions involve the operator $S_b$ which can be the usual Robin transmission operator
or $S_\text{pade}^m$. The reference computation times for solving the Schr\"{o}dinger equation
\eqref{Schrodingereq} on the entire domain are therefore not the same. 
{We denote respectively 
``Robin ref'' and ``$S_{\mathrm{pade}}^{5}$
ref.'' the solution of the Schr\"odinger equation computed on the
whole domain on a single processor for the Robin transmission
condition and for the transmission condition $S_{\mathrm{pade}}^{5}$. }

We see that the new algorithm allows us to reduce the computation time
compared with the classical algorithm. Moreover, the number of
iterations is almost independent of the number of subdomains and the
algorithm is scalable. Finally, the algorithm converges faster with
the transmission condition $S_{\mathrm{pade}}^m$, but takes more
computational time since. Indeed, the application of the transmission condition
$S_{\mathrm{pade}}^m$ is more expensive than the Robin transmission
condition.


\begin{table}[!htbp]
  \footnotesize
  \centering
  \renewcommand{\arraystretch}{1.0}
  \caption{Number of iterations of the first time step and total computation time in seconds, $T=0.5$, $\Delta t=0.01$, $\Delta x=1/128$, $\Delta y=1/8$.}
  \begin{tabular}{|c|c|c|c|c|c|c|}
    \hline
    $N$ & & 2 & 4 & 8 & 16 & 32 \\
    \hline
    & $\mathrm{Robin}^{*}$ & 11 & 11 & 11 & 11 & 11 \\
    \multirow{-2}{*}{Number of iterations}
    & $S_{\mathrm{pade}}^{5}$ & 9 & 9 & 9 & 9 & 10  \\
    \hline
    & Robin ref. & \multicolumn{5}{c|}{16.0} \\                                            
    \cline{3-7}
    & Robin cls. & 63.1 & 32.6 & 17.5 & 11.0 & 5.4 \\
    & Robin new & 30.0 & 9.7 & 4.4 & 2.5 & 1.3 \\
    \cline{2-7}
    & $S_{\mathrm{pade}}^{5}$ ref. & \multicolumn{5}{c|
    }{22.1} \\
    \cline{3-7}
    & $S_{\mathrm{pade}}^{5}$ cls. & 96.1 & 49.8 & 26.7 & 15.0 & 7.8 \\
    \multirow{-6}{*}{Computation time}
    
    & $S_{\mathrm{pade}}^{5}$ new & 38.2 & 14.7 & 6.6 & 3.4 & 1.8 \\
    \hline
  \end{tabular}
  \begin{tablenotes}
  \item *: $p=15$.
  \end{tablenotes}                                       
  \label{Time_clnew_T>dt_128_V0}
\end{table}

%
\begin{remark}
  From the results shown in the following subsection (see Table
  \ref{Niter_m32_V0} and Table \ref{Niter_p32_V0} together), we can
  see clearly the different convergence rates with the two kinds of
  transmission conditions.  
\end{remark}

\begin{figure}[!htbp]
  \centering
  \includegraphics[width=0.46\textwidth]{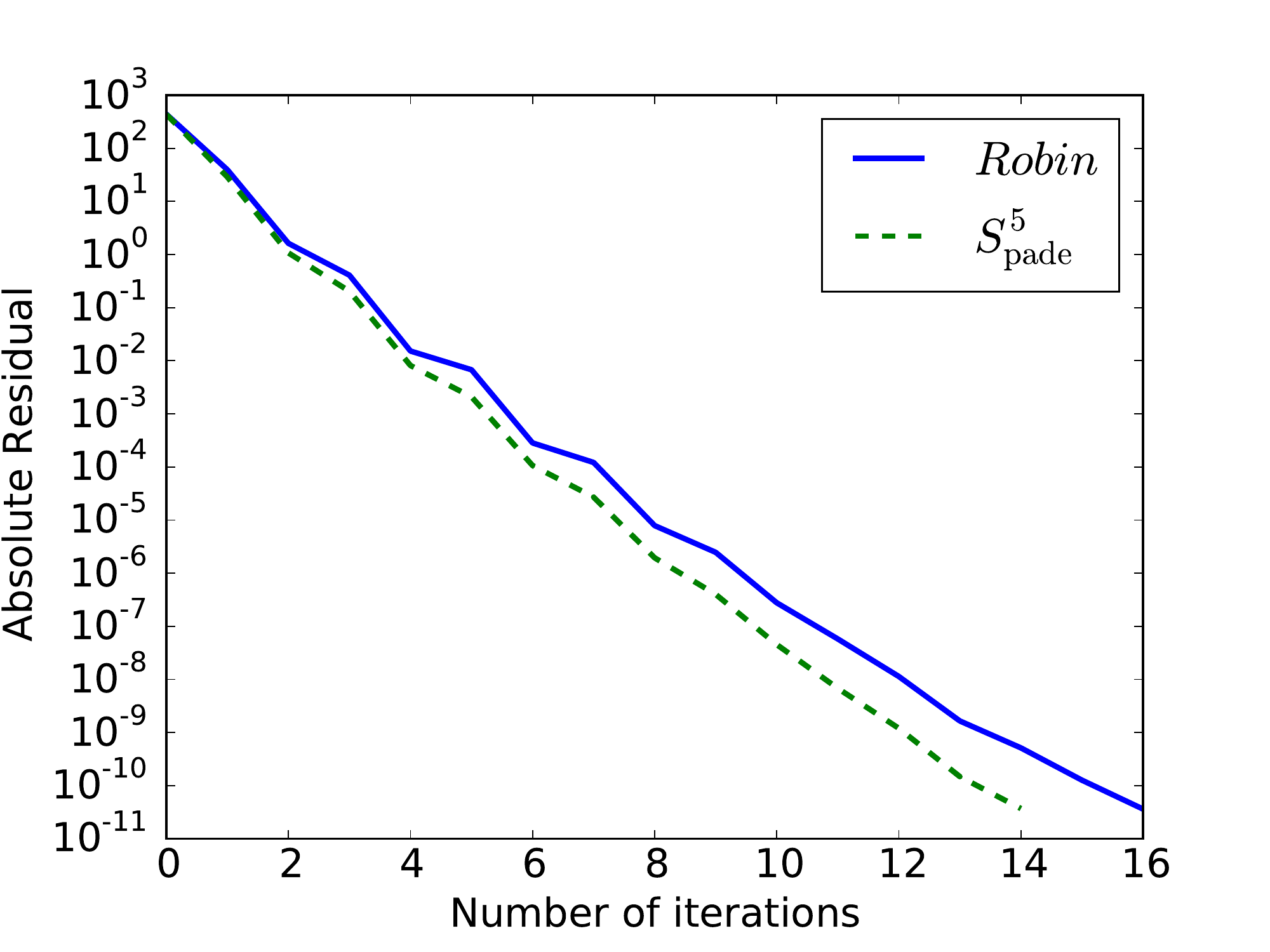}
  \includegraphics[width=0.46\textwidth]{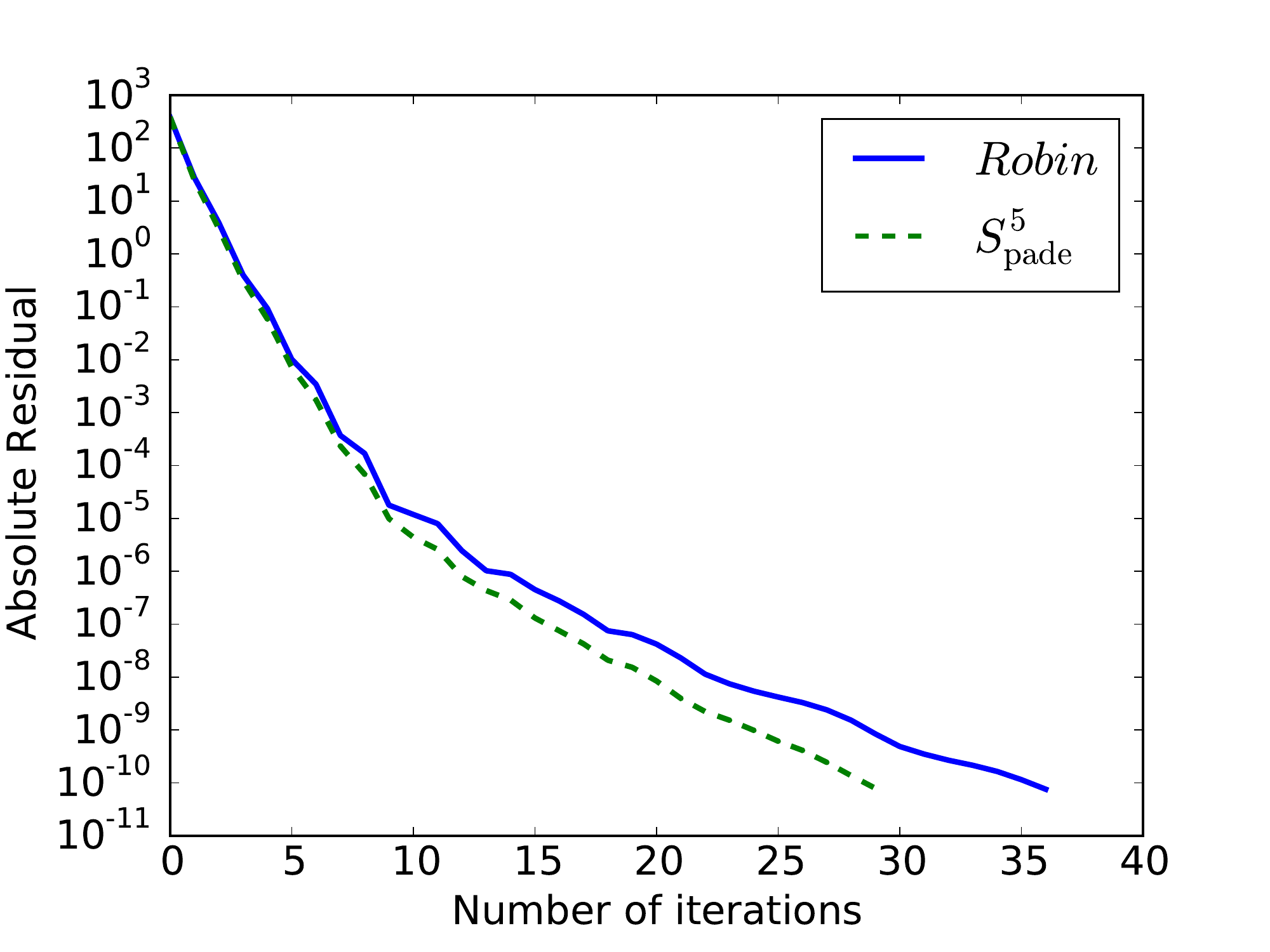}	
  \caption{Convergence history of the first time step for $N=256$ (left), $N=1024$ (right).}
  \label{Hist_256_1024_V0}
\end{figure}

Next we make tests with the mesh $\Delta x = 1/2048$ $\Delta y =
1/64$. The entire domain is divided into $N =128,256,512,1024$
subdomains. We present in Figure \ref{Hist_256_1024_V0} the
convergence history for $N =256,1024$ and in Table
\ref{Time_clnew_T>dt_2048_V0} the numbers of iterations and the total
computation time. We can see that the classical and the new algorithms
are not very scalable since the number of iterations increases with
the number of subdomains. However, the new algorithm takes less
computation time. Besides, the number of iterations required for the
transmission condition $S_{\mathrm{pade}}^m$ is less than the one for the Robin
transmission condition. The computation times are however similar.


\begin{table}[!htbp]
  \footnotesize
  \centering
  \renewcommand{\arraystretch}{1.2}
  \caption{Number of iterations of the first time step and total computation time in seconds, $T=0.5$, $\Delta t=0.01$, $\Delta x=1/2048$, $\Delta y=1/64$.}
  \begin{tabular}{|c|c|c|c|c|c|}
    \hline
    $N$ & & 128 & 256 & 512 & 1024 \\
    \hline
    & $\mathrm{Robin}^{*}$ & 14 & 16 & 22 & 36 \\
    \multirow{-2}{*}{Number of iterations}
    & $S_{\mathrm{pade}}^{5}$ & 12 & 14 & 19 & 29 \\
    \hline
    & Robin cls. & 250.2 & 143.8 & 92.5 & 101.4 \\
    & Robin new & 59.1 & 38.1 & 36.2 & 52.3 \\
    & $S_{\mathrm{pade}}^{5}$ cls. & - & 187.5 & 162.6 & 127.5 \\
    \multirow{-4}{*}{Computation time}	    
    & $S_{\mathrm{pade}}^{5}$ new & - & 42.7 & 36.2 & 45.0 \\
    \hline
  \end{tabular}
  \begin{tablenotes}
  \item -: the memory is insufficient; \hspace{1mm} *: $p=10$.
  \end{tablenotes}
  \label{Time_clnew_T>dt_2048_V0}
\end{table}

\subsubsection{Influence of parameters}
\label{Sec_InfluenceParameter}

We study in this part the influence of parameters: $m$ in the
transmission condition $S_{\mathrm{pade}}^m$ and $p$ in the
transmission condition Robin. The time step is fixed as $\Delta t =
0.01$. The mesh is $\Delta x = 1/128$, $\Delta y = 1/8$. Three
different methods are used to solve the the interface problem: fixed point
method, GMRES method and BiCGStab method. Two initial vectors are
considered: the zero vector and the random vector. In our tests, the
algorithm initialized with the zero vector converges faster than the
one with the random vector. However, in this subsection, our goal is to study
the influence of parameters. As explained in \cite{Gander2008history},
initialization with a zero vector to compute a smooth solution makes
that the error contains only low frequencies and could therefore draw
the wrong conclusions. Thus we consider both of the two initial
vectors, but we have no wish to compare them. 

\vbox{}
\noindent \textbf{The parameter $m$ in the transmission condition $S_{\mathrm{pade}}^m$}.
The numbers of iterations of the first time step with various
$m$ are presented in Table \ref{Niter_m32_V0}. We could see that the
number of iterations is not sensitive to the parameter $m$ if the
GMRES method or the BiCGStab method is used on the interface
problem. However, if the fixed point method is used, increasing the
order of $S_{\mathrm{pade}}^m$ does not ensure better convergence
property. The transmission condition $S_{\mathrm{pade}}^m$ is based on
formal Pad\'{e} approximation of square root operator, this
approximation may deteriorate for large $m$. 

\begin{table}[!htbp] 
  \footnotesize
  \renewcommand{\arraystretch}{1.2}
  \centering
  \caption{Number of iterations for different $m$, $N=32$.}
  \begin{tabular}{|c|c|c|c|c|c|c|}
    \hline
    & \multicolumn{2}{c|}{Fixed point} & \multicolumn{2}{c|}{GMRES} & \multicolumn{2}{c|}{BiCGStab}\\ \hline
    $m$ & Zero & Random & Zero & Random & Zero & Random\\
    \hline
    3 & 34 & 124 & 11 & 30 & 6 & 17 \\
    4 & 26 & 96 & 10 & 28 & 6 & 16 \\
    5 & 21 & 79 & 10 & 27 & 5 & 15 \\
    6 & 18 & 68 & 9 & 26 & 5 & 15 \\
    7 & 17 & 61 & 9 & 25 & 5 & 14 \\
    8 & 18 & 56 & 9 & 25 & 5 & 14 \\
    9 & 19 & 52 & 10 & 25 & 5 & 14 \\
    10 & 21 & 50 & 10 & 25 & 5 & 14\\
    15 & 32 & 46 & 11 & 26 & 6 & 15 \\
    20 & 43 & 50 & 12 & 28 & 6 & 16 \\
    \hline
  \end{tabular}
  \label{Niter_m32_V0}
\end{table}

\vbox{}
\noindent \textbf{The parameter $p$ in the transmission condition Robin}.
We present in Table \ref{Niter_p32_V0} the numbers of iterations with
various $p$ for $N=32$. From the table, we can see that the
algorithm is not sensitive to $p$ if the GMRES method or the BiCGStab
method is used on the interface problem, while it exists an optimal
$p$ (marked with an underline) if the fixed point method is
used. Besides, the Krylov methods could accelerate a lot the
convergence.
%
\begin{table}[!htbp] 
  \footnotesize
  \renewcommand{\arraystretch}{1.2}
  \caption{Number of iterations for different $p$, $N=32$.}                                                        
  \centering                                              
  \begin{tabular}{|c|c|c|c|c|c|c|}
    \hline
    & \multicolumn{2}{c|}{Fixed point} & \multicolumn{2}{c|}{GMRES} & \multicolumn{2}{c|}{BiCGStab}\\ \hline
    $p$ & Zero & Random & Zero & Random & Zero & Random \\
    \hline
    $5$  & 57 & 580 & 11 & 35 & 6 & 21 \\
    $10$ & 34 & 315 & 11 & 32 & 6 & 19 \\
    $15$ & 32 & 239 & 11 & 31 & 6 & 18 \\
    $20$ & 35 & 209 & 11 & 31 & 6 & 19 \\
    $25$ & 40 & 200 & 11 & 32 & 6 & 19 \\
    $30$ & 46 & 199 & 11 & 32 & 6 & 19 \\
    $35$ & 52 & 204 & 12 & 33 & 6 & 19 \\
    $40$ & 59 & 209 & 12 & 33 & 6 & 20 \\
    $45$ & 65 & 222 & 12 & 34 & 6 & 20\\
    \hline
    $15$ & {\underline{32}} & & & & & \\
    $26$ & & {\underline{198}} & & & &\\
    \hline
  \end{tabular}
  \label{Niter_p32_V0}
\end{table}

In conclusion, the difference between the two transmission conditions
is clear if the fixed point method is applied to the interface problem
and the initial vector is the random vector (see Tables
\ref{Niter_m32_V0} and \ref{Niter_p32_V0} together). The number of
iterations for the transmission condition $S_{\mathrm{pade}}^m$ is
less than that for the transmission condition Robin. But the
difference is smaller using the Krylov methods and the zero vector as
the initial vector. From the point of view of computation time, the
zero vector and the GMRES method are a good choice. According to the
tests in the previous subsection, the computation time for the two
transmission conditions are similar. 

\subsection{Case of non-zero potential}
\label{Sec_No0V}

The aim of this section is to compare the classical algorithm and the
preconditioned algorithm in the framework of the DDS method with the
fixed point method used on the interface problem. We first consider
the potential $V=x^2+y^2$. Let us denote by $N_{\mathrm{nopc}}$
(resp. $N_{\mathrm{pc}}$) the number of iterations required to obtain
convergence of the classical algorithm (resp. the preconditioned
algorithm) and $T_{\mathrm{nopc}}$ (resp. $T_{\mathrm{pc}}$) the
computation time of the classical algorithm (resp. the preconditioned
algorithm).  

\begin{figure}[!htbp]
  \centering
  \includegraphics[width=0.46\textwidth]{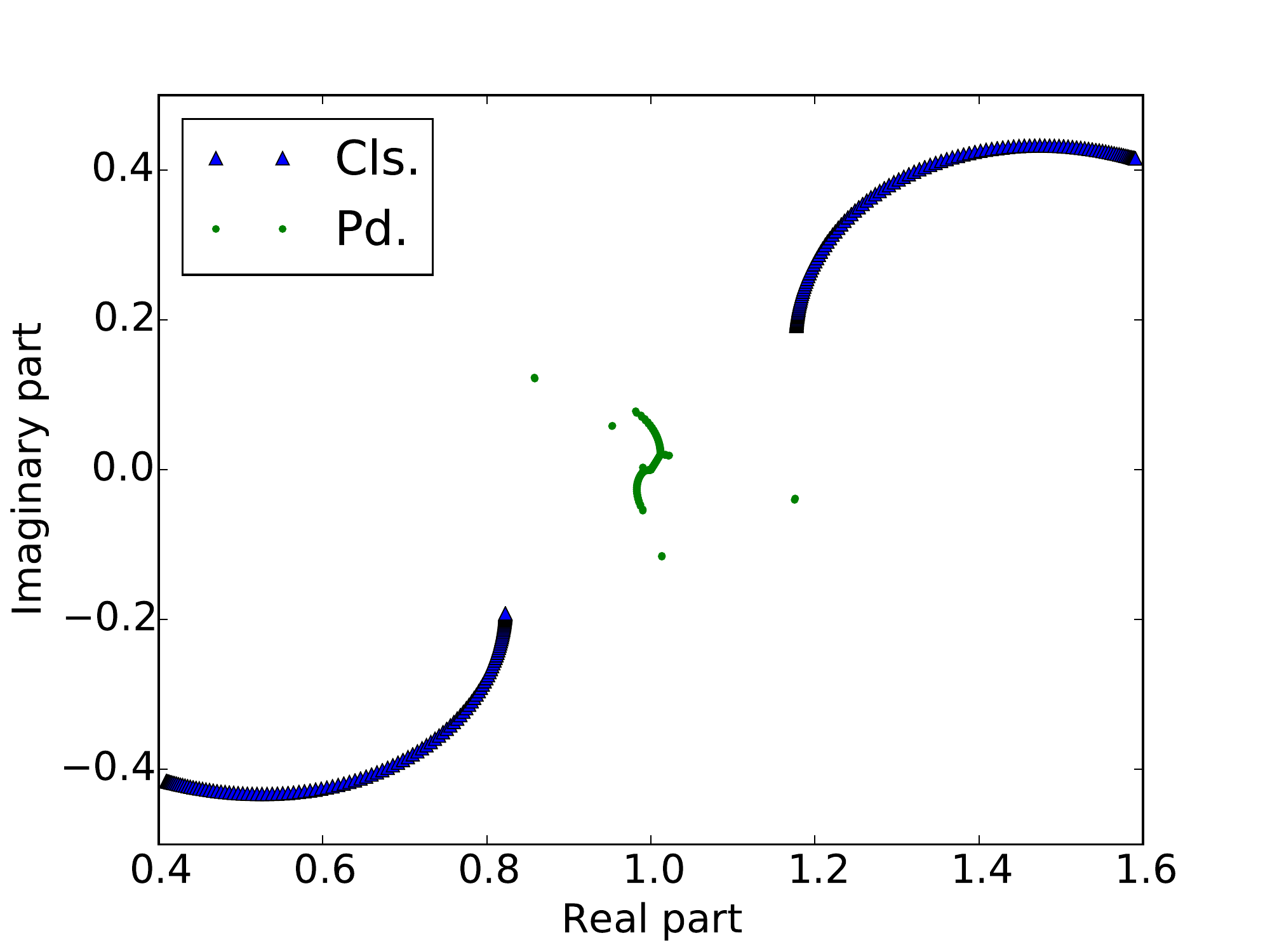}
  \includegraphics[width=0.46\textwidth]{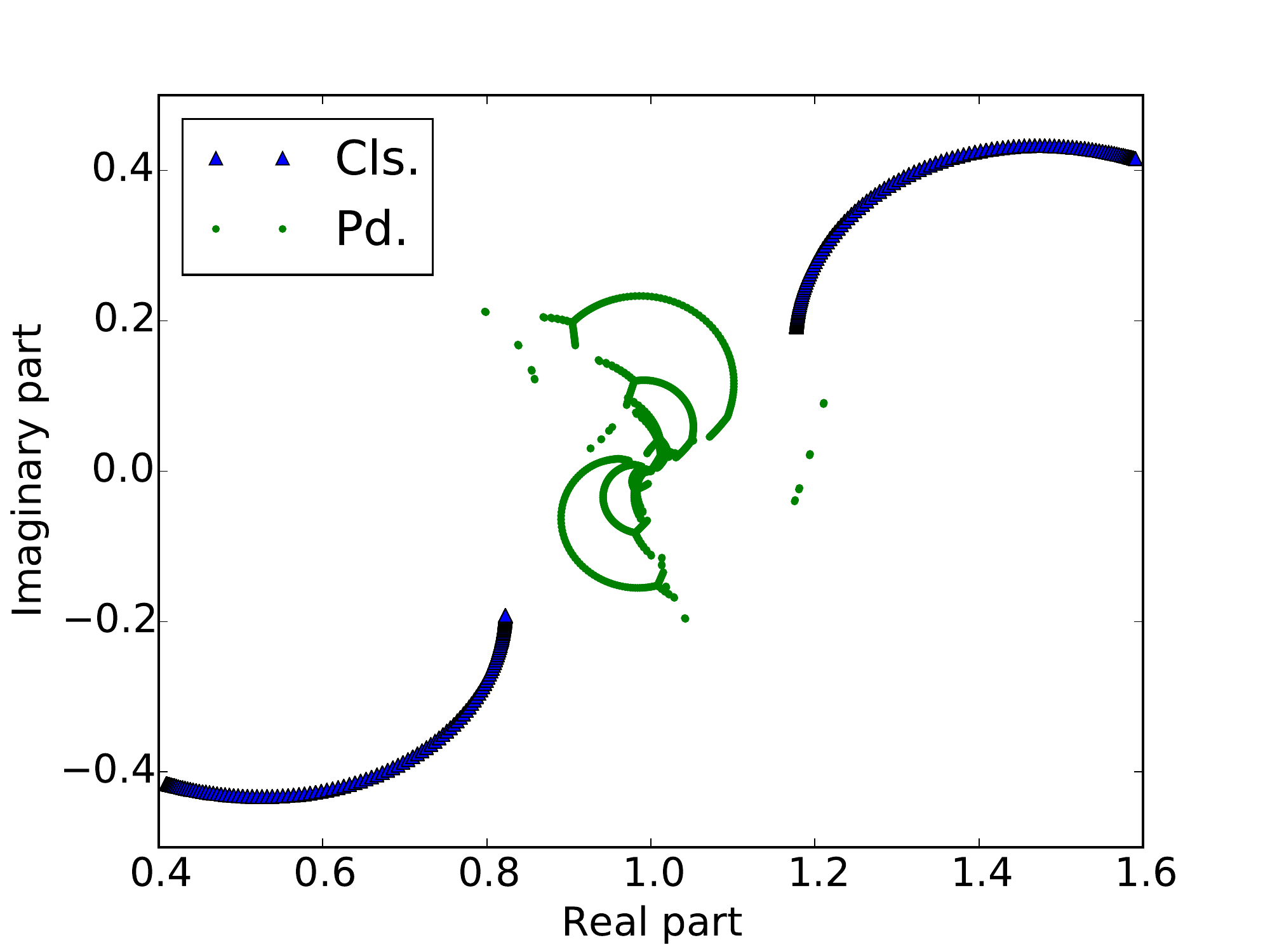}
  \caption{Eigenvalues of the classical (Cls. ) and the preconditioned (Pd.) algorithm, first time step, $\Delta t = 0.01$, $\Delta x = 1/128$, $\Delta y = 1/8$, $N=2$ (left), $N=8$ (right).}
  \label{SpPdN2}
\end{figure}

First, we present in Figure \ref{SpPdN2} the spectral properties of
$(I-\mathcal{L})$ and $P^{-1}(I-\mathcal{L})$ for $N=2,8$
subdomains. We see that the eigenvalues of
  $P^{-1}(I-\mathcal{L})$ are closer to $1+0i$ than those of
  $(I-\mathcal{L})$. The convergence history is therefore better.
The Table \ref{Niter_Time_2_32_L_abc2p_Vxx} shows the number of iterations of
the first time step and the total computation time to realize a
complete simulation. The mesh is $\Delta x=1/128$, $\Delta y=1/8$. As
mentioned before, it is possible to solve the Schr\"{o}dinger equation
on the entire domain with this mesh. The computation time is denoted
by $T^{\mathrm{ref}}$. We can see that all algorithms are robust and
the number of iterations is independent of the number of
subdomains. The classical algorithm and the preconditioned algorithm
are both scalable and the preconditioner allows to reduce the number
of iterations and the total computation time. In addition, the times of
computation $T_{\mathrm{nopc}}$, $T_{\mathrm{pc}}$ are less than the
reference time $T^{\mathrm{ref}}$ if $N$ is large.

\begin{table}[!htbp] 
  \footnotesize
  \centering
  \renewcommand{\arraystretch}{1.2}
  \caption{Number of iterations of the first time step and total computation time in seconds of the classical algorithm and the preconditioned algorithm, $T = 0.5$, $\Delta t = 0.01$, $\Delta x = 1/128$, $\Delta y = 1/8$.}
  \begin{tabular}{|c|c|c|c|c|c|}
    \hline
    $N$ & 2 & 4 & 8 & 16 & 32 \\
    \hline
    $N_{\mathrm{nopc}}$, $m=7$ & 17 & 17 & 17 & 17 & 17 \\
    \hline
    $N_{\mathrm{pc}}$, $m=5$ & 5 & 5 & 5 & 5 & 5 \\
    \hline
    $T^{\mathrm{ref}}$ & \multicolumn{5}{c|}{16.1}\\
    \hline
    $T_{\mathrm{nopc}}$ & 142.7 & 75.3 & 40.1 & 23.9 & 12.1 \\
    \hline
    $T_{\mathrm{pc}}$ & 91.3 & 43.3 & 22.8 & 13.1 & 7.3 \\
    \hline
  \end{tabular}
  \label{Niter_Time_2_32_L_abc2p_Vxx}
\end{table}

\begin{figure}[!htbp]
  \centering
  \includegraphics[width=0.46\textwidth]{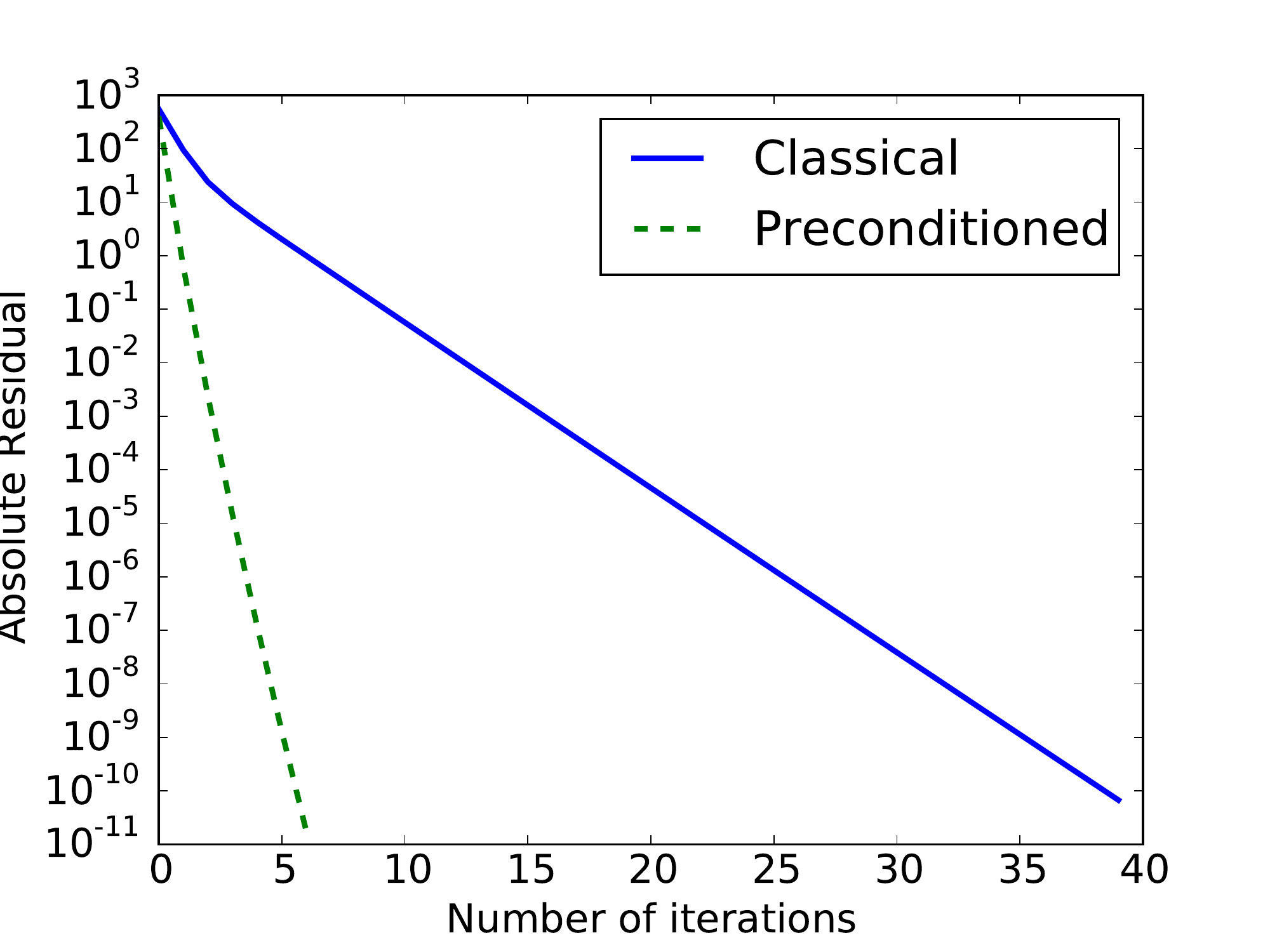}
  \includegraphics[width=0.46\textwidth]{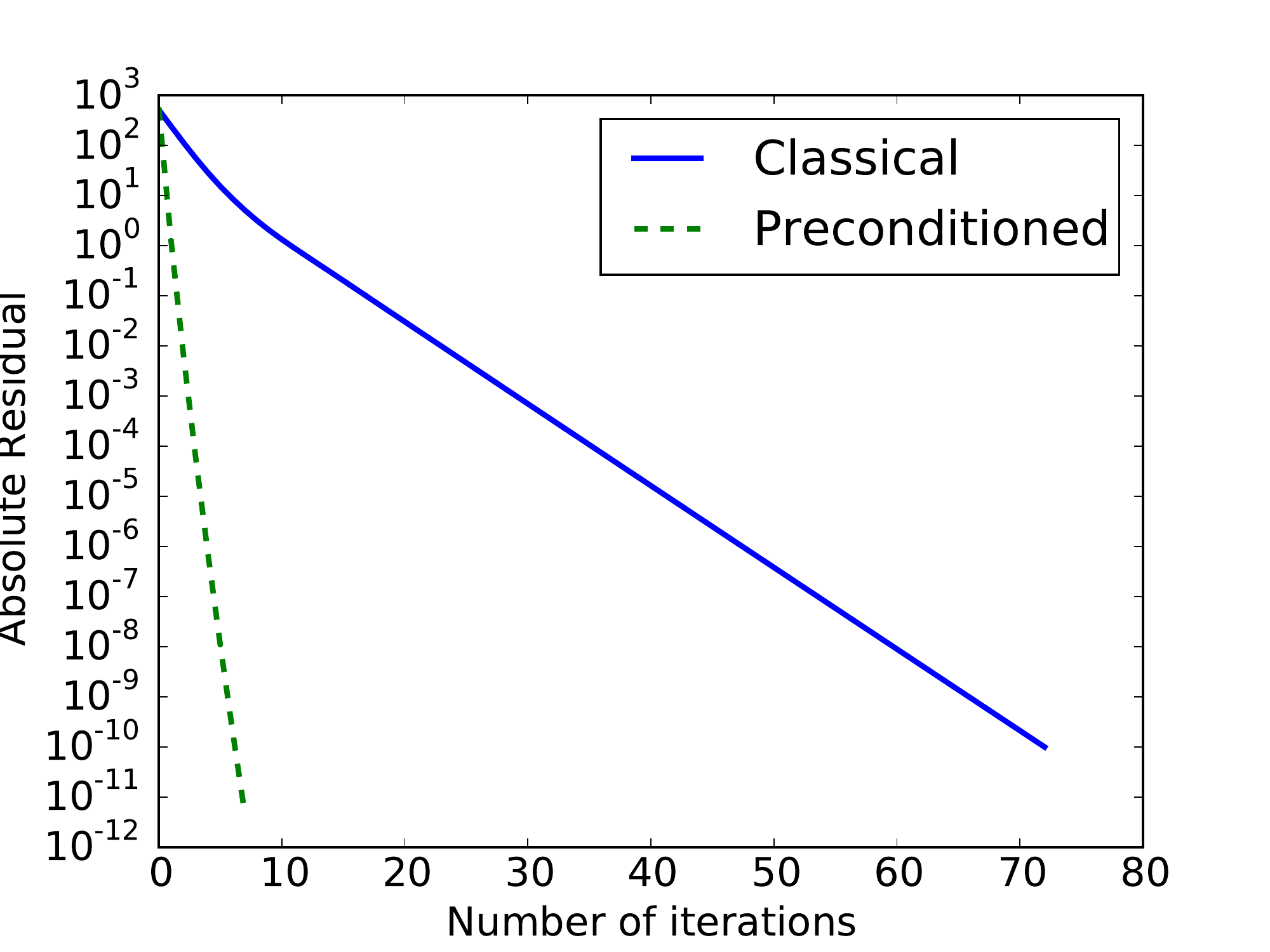}
  \caption{Convergence history of the first time step of the classical
    algorithm and the preconditioned algorithm, $\Delta t = 0.01$,
    $\Delta x = 1/2048$, $\Delta y = 1/64$, $N=256$ (left), $N=512$
    (right).}
  \label{hist_N_256_1024_L_abc2p_Vxx}
\end{figure}

Next, we consider the mesh $\Delta x = 1/2048, \Delta y = 1/64$. The
size of $(I-\mathcal{L})$ is too large to compute all its spectral
values. The convergence history and the computation time are presented
in Figure \ref{hist_N_256_1024_L_abc2p_Vxx} and Table
\ref{Niter_Time_256_1024_L_abc2p_Vxx}. The parameters $m$ used are
also presented in Table \ref{Niter_Time_256_1024_L_abc2p_Vxx}. We can
see that the preconditioned algorithm is more robust and requires much
less number of iterations. However the computation time of the two
algorithms are not quite scalable since for the classical algorithm,
the number of iterations increases with the number of
subdomains. Concerning the preconditioned algorithm, the size of
preconditioner is $(2N-2) \times N_T \times N_y$. This increases with
the number of subdomains $N$. Thus, the application of preconditioner
takes more computation time with bigger $N$. 

%
%

\begin{table}[!htbp] 
  \footnotesize
  \centering
  \renewcommand{\arraystretch}{1.2}
  \caption{Computation time in seconds of the classical algorithm and the preconditioned algorithm, $\Delta t = 0.01, \Delta x = 1/2048, \Delta y = 1/64$.}                                              
  \begin{tabular}{|c|c|c|c|}
    \hline
    $N$ & & 256 & 512 \\
    \hline
    \multirow{2}{*}{Classical algorithm} & $m$ & 8 & 12  \\
    & $T_{\mathrm{nopc}}$ & 417.9 & 344.9 \\
    \hline
    \multirow{2}{*}{Preconditioned algorithm} & $m$ & 5 & 6 \\
    & $T_{\mathrm{pc}}$ & 259.9 & 268.7 \\
    \hline
  \end{tabular}
  \label{Niter_Time_256_1024_L_abc2p_Vxx}
\end{table}

We finish this subsection by some numerical tests for a time dependent
potential $V=5(x^2+y^2)(1+\cos(4\pi t))$. We get similar
conclusion by the results shown in table
\ref{Niter_Time_2_32_L_robin_Vtx}. 
\begin{table}[!htbp] 
  \footnotesize
  \centering
  \renewcommand{\arraystretch}{1.2}
  \caption{Number of iterations of the first time step and total computation time in seconds of the classical algorithm and the preconditioned algorithm, $T = 0.5$, $\Delta t = 0.01$, $\Delta x = 1/128$, $\Delta y = 1/8$.}                                              
  \begin{tabular}{|c|c|c|c|c|c|}
    \hline
    $N$ & 2 & 4 & 8 & 16 & 32 \\
    \hline
    $N_{\mathrm{nopc}}$, $p=-10$ & 31 & 31 & 31 & 31 & 31 \\
    \hline
    $N_{\mathrm{pc}}$, $p=-10$ & 9 &  9 & 9 & 9 & 9 \\
    \hline
    $T^{\mathrm{ref}}$ & \multicolumn{5}{c|}{263.2}\\
    \hline
    $T_{\mathrm{nopc}}$ & 272.6 & 138.5 & 72.4 & 40.5 & 19.5 \\
    \hline
    $T_{\mathrm{pc}}$ & 205.6 & 101.4 & 52.2 & 29.4 & 14.8 \\
    \hline
  \end{tabular}
  \label{Niter_Time_2_32_L_robin_Vtx}
\end{table}

In conclusion, the preconditioner allows to reduce significantly the number of iterations and the computation time.

\section{Conclusion}

We applied the SWR method and the DDS method to the two dimensional
linear Schr\"{o}dinger equation with general potential. We proposed a new algorithm if the
potential is a constant and a preconditioned algorithm for a general
linear potential, which allows to reduce the number of iterations and the
computation time compared with the classical one. According to the
numerical tests, the preconditioned algorithm is not sensitive to the
transmission conditions (Robin, $S_{\mathrm{pade}}^m$) and the
parameters in these conditions.  

\section*{Acknowledgements} We acknowledge Pierre Kestener (Maison de
la Simulation Saclay France) for the discussions about the parallel
implementation. This work was partially supported by the French ANR
grant ANR-12-MONU-0007-02 BECASIM (Mod\`eles Num\'eriques call). The
first author also acknowledges support from the French ANR grant BonD
ANR-13-BS01-0009-01.

\bibliographystyle{abbrv}
\bibliography{Bib2}

\end{document}